\documentclass[lettersize,journal]{IEEEtran}
\bstctlcite{bstctl:etal, bstctl:nodash, bstctl:simpurl}

\IEEEoverridecommandlockouts                              

\usepackage{amsmath,graphicx,amsfonts,amssymb,amsthm,epsfig,mathrsfs,balance,stmaryrd,mathabx}
\usepackage{subcaption,tikz,color,multirow,algorithm,algpseudocode,algorithmicx,tabulary}
\usepackage[noadjust]{cite}
\usepackage{url,framed,bm,dsfont} 
\usepackage{nicefrac}
\usepackage[normalem]{ulem}

\newtheorem{lemma}{Lemma}
\newtheorem{proposition}{Proposition}
\newtheorem{remark}{Remark}

\usepackage{enumitem}  

\newcommand{\tr}{{\mathrm{trace}}}

\newcommand{\differential}{{\rm{d}}}
\newcommand{\diag}{\mathrm{diag}}

\newcommand{\R}{\mathbb{R}}

\renewcommand{\qedsymbol}{\hfill\ensuremath{\blacksquare}}

\DeclareMathOperator*{\proj}{proj}

\newcommand{\red}{\color{red}}

\title{\LARGE \bf
Stochastic Learning of Computational Resource Usage as\\Graph Structured Multimarginal Schrödinger Bridge  
}

\author{Georgiy A. Bondar, Robert Gifford, Linh Thi Xuan Phan, Abhishek Halder,~\IEEEmembership{Senior Member,~IEEE}
\thanks{Georgiy A. Bondar is with the Department of Applied Mathematics, University of California, Santa Cruz, CA 95064, USA, {\texttt{gbondar@ucsc.edu}}.}
\thanks{Robert Gifford and Linh Thi Xuan Phan are with the Department of Computer and Information Science, University of Pennsylvania, Philadelphia, PA 19104, USA, {\texttt{\{rgif,linhphan\}@seas.upenn.edu}}.}
\thanks{Abhishek Halder (corresponding author) is with the Department of Aerospace Engineering, Iowa State University, Ames, IA 50011, USA, {\texttt{ahalder@iastate.edu}}.}         
\thanks{This research was supported by NSF awards 2112755, 2111688 and 1750158.}
}

\begin{document}
\bstctlcite{IEEE_b:BSTcontrol}
\maketitle
\thispagestyle{plain}
\pagestyle{plain}

\begin{abstract}
We propose to learn the time-varying stochastic computational resource usage of software as a graph structured Schr\"{o}dinger bridge problem. In general, learning the computational resource usage from data is challenging because resources such as the number of CPU instructions and the number of last level cache requests are both time-varying and statistically correlated. Our proposed method enables learning the joint time-varying stochasticity in computational resource usage from the measured profile snapshots in a nonparametric manner. The method can be used to predict the most-likely time-varying distribution of computational resource availability at a desired time. We provide detailed algorithms for stochastic learning in both single and multi-core cases, discuss the convergence guarantees, computational complexities, and demonstrate their practical use in two case studies: a single-core nonlinear model predictive controller, and a synthetic multi-core software. 
\end{abstract}




\section{Introduction}\label{sec:Intro}
Compute-intensive software, including control software, often operate in hardware platforms that are different from where their performance were verified. This could be because of hardware up/downgrades in response to evolving project needs, technology progress and such. A natural question, then, is whether there could be a principled way to learn the stochastic dynamical nature of computational resource availability such as processor, memory bandwidth and last level shared cache (LLC). Such learning could then be leveraged to design dynamic scheduling algorithms by predicting, for example, the most likely \emph{joint} computational resource that will be available at a future time. In this work, we propose such a stochastic learning framework for both single and multi-core platforms.

In the real-time systems community, it is well-known \cite{bernat2002wcet} that due to hardware-level stochasticity, multiple executions of the same software on the same hardware with identical initial conditions and parameters, result in different execution times. A common way to account for this is to analyze some lumped variable such as (worst-case or probabilistic) execution time \cite{wilhelm2008worst,davis2011survey,gil2017open} for a given software-hardware combination. In contrast, the ability to directly learn in the joint space of LLC, processor and memory availability, could enable the design of more fine-grained dynamic resource schedulers. This is of particular interest in safety-critical controlled cyber-physical systems that operate in resource-constrained environments because high performance controllers (e.g., MPC) are more compute-intensive and have more pronounced stochasticity than computationally benign controllers (e.g., PID). 

Learning in the joint space, however, is technically challenging because variables such as processor availability, LLC and memory bandwidth are not only statistically correlated, but their correlations also change over time, i.e., they cannot be assumed as (neither strict nor wide-sense) stationary stochastic processes for learning purposes. It is also impractical to fit parametric statistical models such as a Markov chain by gridding the multi-dimensional space of computational resource state since the resource state space is not countable. In this work, we leverage the recent progress in multimarginal Schr\"{o}dinger bridges to show the feasibility of nonparametric joint stochastic learning of computational resource usage from hardware-software profile data.

\subsubsection*{Related work}
Multimarginal Schr\"{o}dinger bridge problems (MSBPs) are entropy regularized variants of the multimarginal optimal mass transport (MOMT) problems \cite{ruschendorf2002n,pass2015multi}. The latter has been applied to learning and inference problems in chemical physics \cite{buttazzo2012optimal,cotar2013density}, team matching \cite{carlier2015numerical}, fluid dynamics \cite{benamou2019generalized} and risk management \cite{ennaji2024robust}. For recent applications of MSBPs to learning and tracking problems, see e.g., \cite{elvander2020multi,haasler2021multimarginal,noble2024tree}. The MSBPs can also be seen as generalized variants of the classical (a.k.a. bimarginal) Schr\"{o}dinger bridge problems (SBPs) which we will explain in Sec. \ref{subsec:bimarginalSBP}. In the control literature, there are growing works on the stochastic optimal control interpretations \cite{chen2021stochastic,teter2023contraction}, generalizations \cite{caluya2020finite,caluya2021reflected,caluya2021wasserstein} and applications \cite{nodozi2022schrodinger,haddad2022density,nodozi2023physics} of the SBPs.

Prior work has applied machine learning to performance prediction and workload modeling, especially in data center and cloud environments; see e.g., the surveys in~\cite{amiri2017survey, Saxena23survey} and references therein. Unlike our work, these solutions focus on predicting {\em coarse-grained} characteristics -- such as request rates, CPU utilization, memory and disk I/O usages, bandwidth, or energy consumption -- and they do not consider the interdependent relationship among shared resources. In real-time and embedded systems, learning techniques have been used for estimating timing behaviors~\cite{Kumar2022,Jaekwon2023,muts2022multiobjective}; however, existing work focuses on the worst-case or probabilistic execution time of the entire software, rather than the 
dynamic 
resource usage patterns during an execution of the software. To our best knowledge, our work is the first to use learning for predicting the stochastic dynamical run-time behavior of computational resource usage that jointly considers interdependent resources on multi-core platforms (such as CPU, cache and memory bandwidth) as well as multi-threaded applications.

\subsubsection*{Contributions} 
This work builds on our preliminary work \cite{bondar2023path} but significantly extends the same by considering multi-core resources and more general graph structures than path tree. The specific contributions are threefold:

\begin{itemize}
    \item Mathematical formulations (Sec. \ref{sec:StocLearningOfHardwareResource}-\ref{sec:GraphStructuredMSBP}) for the stochastic learning of computational resource usage as graph-structured MSBPs where the graph structures arise from single or multi-core computational platforms. These formulations are motivated by maximum likelihood interpretations in the path space.

    \item Numerical algorithms (Sec. \ref{sec:Algorithms}) for solving the proposed graph-stuctured MSBPs. We show that unlike generic graph structures, the MSBPs induced by single or multi-core computational resource usage lead to structured tensor optimization problems that can be leveraged to reduce the learning algorithms' computational complexity from exponential to linear in number of observations.

    \item Numerical experiments (Sec. \ref{sec:Experiements}) to illustrate the proposed learning framework. This includes a single-core MPC benchmark, and a synthetic multi-core software benchmark.
\end{itemize}



\section{Preliminaries}\label{sec:prelim}
In this Section, we fix notations and background ideas that will be used subsequently.

\noindent\textbf{Sets.}
For any natural number $\nu$, we use the finite set notation $\llbracket \nu\rrbracket := \{1,2,\hdots,\nu\}$. We denote the cardinality of set $\mathcal{S}$ as $\vert\mathcal{S}\vert$. We will use the following.
\begin{lemma}\cite[Thm. 19.2]{warner1990modern}\label{LemmaCardnilaitySetMinus}
Given finite sets $\mathcal{S}_{1}\subseteq\mathcal{S}_{2}$ with $\vert\mathcal{S}_1\vert=\nu_1$ and $\vert\mathcal{S}_2\vert=\nu_2$, we have $\vert\mathcal{S}_{2}\setminus\mathcal{S}_{1}\vert = \nu_2 - \nu_1$.
\end{lemma}

\noindent\textbf{Vectors, matrices, tensors.}
We use unboldfaced (resp. boldfaced) small letters
to denote scalars (resp. vectors). Unboldfaced capital letters denote matrices and
bold capital letters denote tensors of order three or
more. We occasionally make exceptions for standard notations such as \eqref{DefKLdivergence}-\eqref{DefWasserstein}.

We use square braces to denote the components. For example, $\left[\bm{X}_{i_{1},\hdots,i_{r}}\right]$ denotes the $(i_1,\hdots,i_{r})$th component of an order $r$ tensor $\bm{X}$, where $(i_1,\hdots,i_{r})\in\mathbb{N}^{r}$. We use the $r$ fold tensor product space notation $\left(\mathbb{R}^{d}\right)^{\otimes r} := \underbrace{\mathbb{R}^{d} \otimes \hdots \otimes\mathbb{R}^{d}}_{r\;\text{times}}$. 

For two given tensors $\bm{X},\bm{Y}$ of order $r$, their \emph{Hilbert-Schmidt inner product} is
\begin{align}
\langle\bm{X},\bm{Y}\rangle := \sum_{i_1,\hdots,i_r}\left[\bm{X}_{i_{1},\hdots,i_{r}}\right]\left[\bm{Y}_{i_{1},\hdots,i_{r}}\right].
\label{HilbertSchmidtInnerProdTensors}
\end{align}
The operators $\exp(\cdot)$, $\log(\cdot)$, $\odot$ and $\oslash$ are all understood \emph{elementwise}, i.e., denote elementwise exponential, logarithm, multiplication and division, respectively. Vector or matrix transposition is denoted by the superscript $^{\top}$. We use ${\rm{diag}}(\bm{v})$ to denote a diagonal matrix with entries of vector $\bm{v}$ along its main diagonal, and $\bm{1}$ to denote all ones column vector of suitable length.

\noindent\textbf{Probability.} A \emph{probability measure} $\mu$ over some Polish space $\mathcal{X}$ satisfies $\int_{\mathcal{X}}\differential\mu=1$. For a pair of probability measures $\mu,\nu$ defined over two Polish spaces $\mathcal{X},\mathcal{Y}$ respectively, their \emph{product measure} is $\mu\otimes\nu$, and $\int_{\mathcal{X}\times\mathcal{Y}}\differential\left(\mu\otimes\nu\right)=1$.

The \emph{entropy} of a probability measure $\mu$ is $-\int\log\mu\differential\mu$. The \emph{relative entropy} or \emph{Kullback-Leibler divergence} $D_{\rm{KL}}(\cdot\Vert\cdot)$ between two probability measures $\mu$ and $\nu$ is
\begin{align}
D_{\rm{KL}}(\mu\Vert\nu):=\begin{cases}
\int \log\frac{\differential\mu}{\differential\nu}\differential\mu & \text{if}\quad\mu \ll \nu,\\
+\infty &\text{otherwise,}
\end{cases}
\label{DefKLdivergence}
\end{align}
where $\frac{\differential\mu}{\differential\nu}$ denotes the Radon-Nikodym derivative, and $\mu \ll \nu$ is a shorthand for ``$\mu$ is absolutely continuous w.r.t. $\nu$".

The Wasserstein distance $W$ between two probability measures $\mu,\nu$, supported respectively on $\mathcal{X},\mathcal{Y}\subseteq\mathbb{R}^{d}$, is 
\begin{align}
W\left(\mu,\nu\right):=\left(\underset{\pi\in\Pi(\mu,\nu)}{\inf}\int_{\mathcal{X}\times\mathcal{Y}}\|\bm{x}-\bm{y}\|_2^2\:\differential\pi(\bm{x},\bm{y})\right)^{\frac{1}{2}},
\label{DefWasserstein}    
\end{align}
where the infimum in \eqref{DefWasserstein} is over all joint couplings of $\mu,\nu$, i.e., $\Pi(\mu,\nu):=\{\text{probability measures}\,\pi \mid \int_{\mathcal{Y}}\differential\pi = \mu, \int_{\mathcal{X}}\differential\pi = \nu\}$. Unlike the Kullback-Leibler divergence \eqref{DefKLdivergence}, the Wasserstein distance \eqref{DefWasserstein} is a metric on the space of probability measures. 

\noindent\textbf{Hilbert projective metric.} The Hilbert projective metric \cite{birkhoff1957extensions,bushell1973hilbert,kohlberg1982contraction} $d_{\rm{H}}\left(\bm{u},\bm{v}\right)$ between $\bm{u},\bm{v}\in\mathbb{R}^{n}_{>0}$ (positive orthant) is
\begin{align}
d_{\rm{H}}\left(\bm{u},\bm{v}\right) = \log\left(\dfrac{\max_{i=1,\hdots,n}u_i/v_i}{\min_{i=1,\hdots,n}u_i/v_i}\right).
\label{HilbertMetricPosOrthant}
\end{align}
The convergence plots w.r.t. Hilbert metric for our numerical experiments will be given in Sec. \ref{sec:Experiements}.

The formula \eqref{HilbertMetricPosOrthant} is a special case of the general definition of Hilbert metric $d_{\rm{H}}$ between any two elements of a pointed\footnote{A cone $\mathcal{K}$ is pointed if $\mathcal{K}\cap -\mathcal{K}=\{0\}$.} convex cone $\mathcal{K}$ in a real vector space:
\begin{align}
d_{\rm{H}}(u,v) = \log\left(\dfrac{\inf\{\lambda \geq 0 \mid \lambda v \geq u\}}{\sup\{\lambda \geq 0 \mid u \geq \lambda v\}}\right)\quad\forall u,v\in\mathcal{K}.
\label{HilbertMetricGeneric}    
\end{align}
Note that \eqref{HilbertMetricGeneric} reduces to \eqref{HilbertMetricPosOrthant} for $\mathcal{K}\equiv\mathbb{R}^{n}_{>0}$.


\section{Stochastic Learning of Computational Resource}\label{sec:StocLearningOfHardwareResource}
We collect distributional data for the \emph{computational resource availability state} thought of as a continuous-time stochastic process over time $\tau\in[0,t]$. In what follows, we will make these ideas precise. For now, let us begin by assuming that the data are collected at $s\in\mathbb{N}, s\geq 2$ snapshots
\begin{align}
\tau_1 \equiv 0<\tau_2<\ldots<\tau_{s-1}<\tau_s \equiv t.
\label{SnapshotTimes}
\end{align}
Define the \emph{snapshot index set} $\llbracket s\rrbracket :=\{1,2,\hdots,s\}$. 

We consider multi-core computing hardware with $J\in\mathbb{N}$ cores, and define the \emph{core index set}
$\llbracket J\rrbracket := \{1,2,\hdots,J\}$. For $j\in\llbracket J\rrbracket$ and for a given time horizon of interest $[0,t]$, we think of the \emph{computational resource state vector} $\bm{\xi}^{j}(\tau)$, $0\leq \tau\leq t$, as a continuous time $\mathbb{R}^{d}$-valued stochastic process. 

For each $j\in\llbracket J\rrbracket$, a (random) vectorial sample path $\bm{\xi}^{j}(\tau)$ is referred to as a \emph{profile}. For the single core ($J=1$) case, we drop the core index superscript $j$ and simply denote the sample path or profile as $\bm{\xi}(\tau)$.

As a concrete example, consider the case $\bm{\xi}^{j}\in\mathbb{R}^{3}$ $\forall j \in\llbracket J\rrbracket$, with the components of $\bm{\xi}^{j}$ as
\begin{align}
\begin{pmatrix}
\xi_1^{j}\\
\xi_2^{j}\\
\xi_3^{j}
\end{pmatrix}=
\begin{pmatrix}
\text{instructions retired}\\
\text{LLC requests}\\
\text{LLC misses}
\end{pmatrix}\:\forall j \in\llbracket J\rrbracket.
\label{defFeatureVec} 
\end{align}
In this example, the three elements of $\bm{\xi}^{j}$ denote the number of CPU instructions, the number of LLC requests, and the number of LLC misses in the last time unit (e.g., in the last 10 ms if the profiling frequency is 100 Hz), respectively, in the core $j\in\llbracket J\rrbracket$. We will use this specific $\bm{\xi}^{j}\in\mathbb{R}^{3}$ in our numerical experiments (both single and multi-core). However, we will develop the proposed stochastic learning framework for generic $\bm{\xi}^{j}\in\mathbb{R}^{d}$. This will allow generalizability in the sense that the proposed method can be adapted to an application at hand with a custom definition of the $d\in\mathbb{N}$ components of the stochastic state $\bm{\xi}^{j}$.

\begin{figure*}[h]
    \centering    \includegraphics[width=0.9\linewidth]{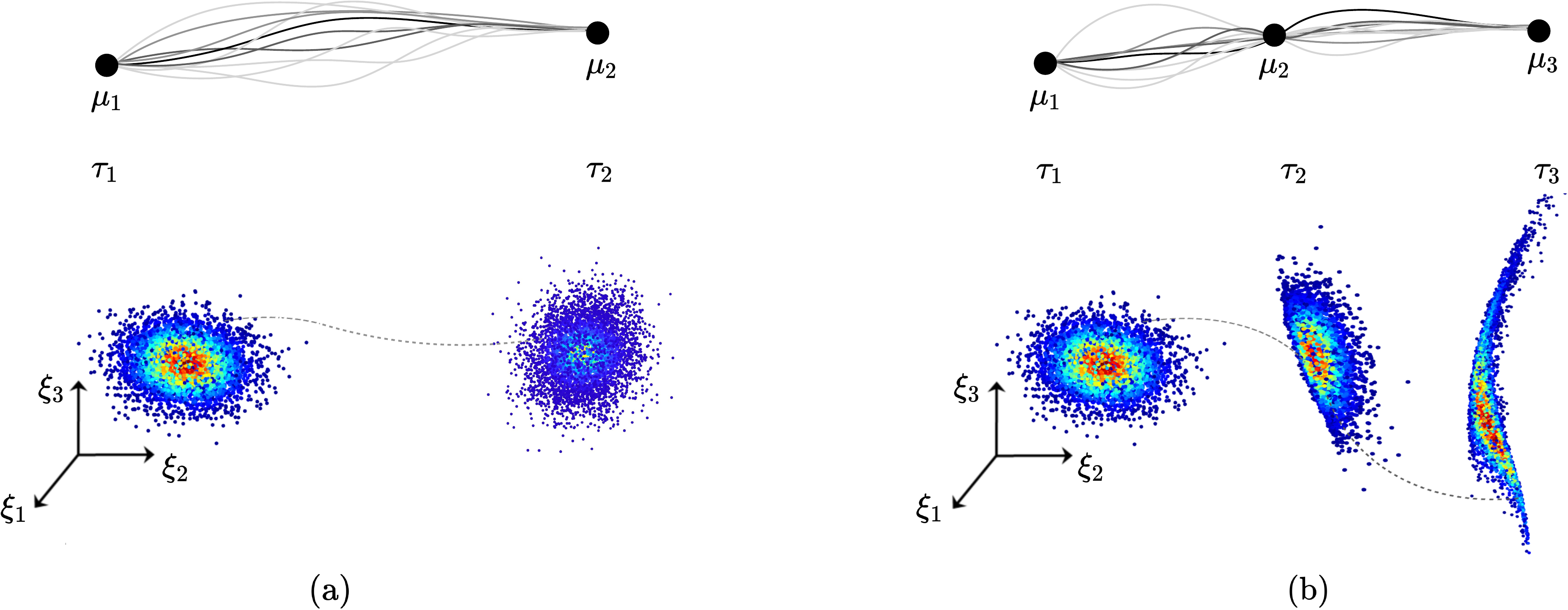}    \caption{{\small{(a) The classical (bimarginal) SBP computes the most likely measure-valued curve connecting two given probability measures $\mu_1,\mu_2$ at times $\tau_1,\tau_2$ respectively. (b) The MSBP computes the most likely measure-valued curve connecting multiple, here three given probability measures $\mu_1,\mu_2, \mu_3$ at times $\tau_1,\tau_2,\tau_3$ respectively. In both subfigures, the feasible measured-valued curves are shown in the top with darker (resp. lighter) hues for higher (resp. lower) probability. The given measures are shown in the bottom as colored scatter plots (red = high probability, blue = low probability) in the ground coordinates $(\xi_1,\xi_2,\xi_3)^{\top}\in\mathbb{R}^3$.}}}
\vspace*{-0.15in}
\label{SBPconcept}
\end{figure*}

The distributional (i.e., measure-valued) observations collected at instances \eqref{SnapshotTimes} comprise a sequence of joint state probability distributions or measures $\{\mu_{\sigma}^{j}\}_{(j,\sigma)\in \llbracket J\rrbracket\times \llbracket s \rrbracket}$, i.e.,
\begin{align}
\bm{\xi}^{j}\left(\tau_\sigma\right) \sim \mu_\sigma^{j},\int \mathrm{d} \mu_\sigma^{j}\left(\boldsymbol{\xi}^{j}\left(\tau_\sigma\right)\right)=1 \quad \forall (j,\sigma) \in \llbracket J\rrbracket\times \llbracket s \rrbracket.
\label{DefObservationSeq}
\end{align}
In this work, we are interested in learning the hardware-level stochasticity, i.e., the stochasticity in the state $\bm{\xi}^{j}$ arises from dynamic resource variability. In other words, repeated execution of the same (e.g., control) software with the same initial condition and same parameters result in different profiles $\bm{\xi}^{j}(\tau)$, $(j,\sigma) \in \llbracket J\rrbracket\times \llbracket s \rrbracket$, $0\leq \tau\leq t$. See Fig. \ref{FigAllMeasuredStatesForOneContext} and Fig. \ref{FigAllMeasuredStatesCanneal} as exemplar single core and multi-core profiles, respectively. 

Intuitively, when the hardware-level stochasticity is negligible, then the distributions $\mu_{\sigma}^{j}$ will be approximately Dirac deltas supported on the graph of a single path. When the stochastic variability is significant, $\mu_{\sigma}^{j}$ will have significant dispersion.

In practice, the probability measures $\{\mu_{\sigma}^{j}\}_{(j,\sigma) \in \llbracket J\rrbracket\times \llbracket s \rrbracket}$ are only available empirically from a fixed, say $n\in\mathbb{N}$ profiles. Let the sample or profile index $i\in\llbracket n\rrbracket$. We accordingly set
\begin{align}
\mu_\sigma^{j}:=\frac{1}{n} \sum_{i=1}^n \delta\left(\boldsymbol{\xi}^{j}-\boldsymbol{\xi}^{i,j}\left(\tau_\sigma\right)\right) \: \forall (j,\sigma) \in \llbracket J\rrbracket\times \llbracket s \rrbracket
\label{DefEmpiricalMeasure}    
\end{align}
where $\delta\left(\boldsymbol{\xi}^{j}-\boldsymbol{\xi}^{i,j}\left(\tau_\sigma\right)\right)$ denotes the Dirac delta at the $i$\textsuperscript{th} sample location $\boldsymbol{\xi}^{i,j}\left(\tau_\sigma\right)\in\mathbb{R}^{d}$ for a fixed index pair $(j,\sigma) \in \llbracket J\rrbracket\times \llbracket s \rrbracket$. For any fixed pair $(j,\sigma) \in \llbracket J\rrbracket\times \llbracket s \rrbracket$, the finite set $\{\boldsymbol{\xi}^{i,j}\left(\tau_\sigma\right)\}_{i=1}^{n}$ is a scattered point cloud.

With the basic notations in place, the informal statement for our stochastic learning problem is as follows.\\
\noindent\textbf{Most likely distributional learning problem (informal).} Given distributional snapshots \eqref{DefObservationSeq}, predict the most likely distribution of the computational resource state
\begin{align}
\bm{\xi}^{j}\left({\tau}\right)\sim\mu_\tau^{j}, \;j\in\llbracket J\rrbracket, \;\text{for any}\:\tau\in[0,t].
\label{ProblemStatementInformal} \end{align}
In the next section, we will formalize this problem statement and discuss the related information graph structures.


\section{MSBP and Graph Structures}\label{sec:GraphStructuredMSBP}
To motivate the mathematical formulation, we start by outlining the classical (bimarginal) SBP in Sec. \ref{subsec:bimarginalSBP}. Then in Sec. \ref{subsec:MSBPformulation}, we make the informal statement mentioned at the end of the previous Section rigorous using the framework of large deviation principle \cite{dembo2009large}. The resulting MSBP formulation takes the form of a maximum likelihood problem in the space of probability measure-valued curves, generalizing a similar formulation for the classical SBP. 

We then consider specific cases of this MSBP formulation which result from the information graph structures induced by the stochastic profiles in single and multi-core computing hardware. In particular, Sec. \ref{subsec:PathStructuredMSBP} details MSBP over a path tree which is the information graph structure that arises in the single core profiling. For the multi-core case, we discuss two different information graph structures in Sec. \ref{subsec:BaryStructuredMSBP} and \ref{subsec:SeriesParallelGraphStructuredMSBP}.      

\subsection{Classical SBP and its Maximum Likelihood Iterpretation}\label{subsec:bimarginalSBP}
In 1931-32, Erwin Schr\"{o}dinger formulated what is now called the classical SBP, in two works: one written in German \cite{schrodinger1931umkehrung} and another in French \cite{schrodinger1932theorie}. For a relatively recent survey, see \cite{leonard2014survey}. The classical SBP corresponds to the \emph{bimarginal} (i.e., $s=2$) \emph{maximum likelihood problem}: its solution finds the \emph{most likely measure-valued curve} $\mu_{\tau}$ where $\tau\in[0,t]$ connecting the given endpoint measures $\mu_1,\mu_2$ at $\tau_1 = 0$ and $\tau_2=t$ with respective supports over subsets of $\mathbb{R}^{d}$. See Fig. \ref{SBPconcept}(a).

The classical SBP formulation proceeds as follows. Let $\mathcal{X}_{1},\mathcal{X}_{2}\subseteq\mathbb{R}^{d}$ be the supports of the given endpoint measures $\mu_1,\mu_2$, respectively. Letting $\bm{\mathcal{X}}:=\mathcal{X}_1\times\mathcal{X}_2 \subseteq \mathbb{R}^{d}\times\mathbb{R}^{d}$, we define a symmetric ground cost $C:\bm{\mathcal{X}}\mapsto \mathbb{R}_{\geq 0}$, i.e., $C(\bm{\xi}(\tau_1),\bm{\xi}(\tau_2))$ is a distance between the random vectors $\bm{\xi}(\tau_1),\bm{\xi}(\tau_2)$. Let  $\mathcal{M}\left(\bm{\mathcal{X}}\right)$ denote the manifold of joint probability measures on the product space $\bm{\mathcal{X}}$. For fixed but not necessarily small $\varepsilon>0$, the classical SBP is an infinite dimensional convex problem:
\begin{subequations}
\begin{align}
&\underset{\bm{M}\in \mathcal{M}\left(\bm{\mathcal{X}}\right)}{\min} \!\int_{\bm{\mathcal{X}}}\!\big\{C\left(\bm{\xi}(\tau_{1}),\bm{\xi}(\tau_{2})\right)+ \varepsilon\log\bm{M}\!\left(\bm{\xi}(\tau_{1}),\bm{\xi}(\tau_{2})\right)\!\!\big\}\nonumber\\
&\qquad\qquad\qquad\qquad\qquad\qquad\qquad\differential\bm{M}\!\left(\bm{\xi}(\tau_{1}),\bm{\xi}(\tau_{2})\right)\label{SBPobj}\\
&\text{subject to}\int_{\mathcal{X}_{2}}\differential\bm{M}\!\left(\bm{\xi}(\tau_{1}),\bm{\xi}(\tau_{2})\right) = \mu_{1},
\label{SBPconstr1}\\
&\qquad\qquad\int_{\mathcal{X}_{1}}\differential\bm{M}\!\left(\bm{\xi}(\tau_{1}),\bm{\xi}(\tau_{2})\right) = \mu_{2}.
\label{SBPconstr2}
\end{align}
\label{SBP}
\end{subequations}
In words, \eqref{SBP} seeks to compute a joint probability measure that minimizes the entropy-regularized transportation cost \eqref{SBPobj} subject to prescribed marginal constraints \eqref{SBPconstr1}-\eqref{SBPconstr2}.

Notice that $\mathcal{M}(\bm{\mathcal{X}})$ is a convex set. The objective \eqref{SBPobj} is strictly convex in $\bm{M}$, thanks to the $\varepsilon$-regularized negative entropy term $\int_{\bm{\mathcal{X}}}\varepsilon\log\bm{M}\:\differential\bm{M}$. The constraints \eqref{SBPconstr1}-\eqref{SBPconstr2} are linear.

To clarify the maximum likelihood interpretation of \eqref{SBP}, let $\mathcal{C}\left([\tau_1,\tau_2],\mathbb{R}^{d}\right)$ denote the collection of continuous functions on the time interval $[\tau_1,\tau_2]$ taking values in $\mathbb{R}^{d}$. Let $\Pi(\mu_1,\mu_2)$ be the collection of all path measures on $\mathcal{C}\left([\tau_1,\tau_2],\mathbb{R}^{d}\right)$ with time $\tau_1$ marginal $\mu_1$, and time $\tau_2$ marginal $\mu_2$. Given a symmetric ground cost (e.g., Euclidean distance) $C:\mathcal{X}_1\times\mathcal{X}_2\mapsto\mathbb{R}_{\geq 0}$, let 
\begin{align}
K(\cdot,\cdot):=\exp\left(-\dfrac{C(\cdot,\cdot)}{\varepsilon}\right),
\label{defkBi}
\end{align} 
and consider the \emph{bimarginal Gibbs kernel} 
\begin{align}
K\left(\bm{\xi}(\tau_1),\bm{\xi}(\tau_2)\right)\mu_1\otimes\mu_2.
\label{DefGibbsKernelBi} 
\end{align}
Recall that the minimizer of problem \eqref{SBP} is an optimal coupling $\bm{M}_{\rm{opt}}\in\mathcal{M}\left(\bm{\mathcal{X}}\right)$. Proposition \ref{prop:SBPmaximumlikelihood} next formalizes why the solution of \eqref{SBP} corresponds to the most likely measure-valued path (Fig. \ref{SBPconcept}(a)) consistent with the observed measure-valued snapshots $\mu_1,\mu_2$. Its proof uses Sanov's theorem \cite{sanov1958probability}; for details we refer the readers to the references below.

\begin{proposition}\label{prop:SBPmaximumlikelihood}(\cite[Sec. II]{follmer1988random},\cite[Sec. 2.1]{pavon2021data})
The optimal coupling $\bm{M}_{\rm{opt}}\in\mathcal{M}\left(\bm{\mathcal{X}}\right)$ in \eqref{SBP} corresponds to the optimal measure-valued path $\pi_{\rm{opt}}\in\Pi(\mu_1,\mu_2)$ solving the (scaled) relative entropy minimization problem:
\begin{align}
\underset{\pi\in\Pi(\mu_1,\mu_2)}{\min}\varepsilon D_{\rm{KL}}\left(\pi\Vert K\left(\bm{\xi}(\tau_1),\bm{\xi}(\tau_2)\right)\mu_1\otimes\mu_2\right).
\label{BimarginalKL}
\end{align}
\end{proposition}

\begin{remark} (\textbf{Existence-uniqueness of the solution for \eqref{BimarginalKL}})\label{Remark:ExistenceUniquenessBimarginal}
Under the stated assumptions on the ground cost $C$, the existence of minimizer for \eqref{BimarginalKL} is guaranteed \cite{csiszar1975divergence,borwein1994entropy}. The uniqueness of the minimizer follows from strict convexity of the map $\pi \mapsto D_{\rm{KL}}(\pi\Vert\nu)$ for fixed $\nu$.
\end{remark}
\begin{remark}
Intuitively, Proposition \ref{prop:SBPmaximumlikelihood} links the solution of a static optimization problem \eqref{SBP} with that of a dynamic optimization problem \eqref{BimarginalKL} in the sense that the minimization in \eqref{BimarginalKL} is performed over all continuous measure-valued curves connecting the endpoints $\mu_1,\mu_2$. 
\end{remark}
Our setting in \eqref{ProblemStatementInformal} requires generalizing these ideas for the multimarginal ($s\geq 2$) case discussed next.

\subsection{MSBP Formulation}\label{subsec:MSBPformulation}
In the single core $(J=1)$ case, Fig. \ref{SBPconcept}(b) illustrates how the MSBP of our interest generalizes the classical a.k.a. bimarginal SBP in Fig. \ref{SBPconcept}(a). 

In general, for any $J\in\mathbb{N}$, we start by defining 
\begin{align}
\mathcal{X}_{\sigma}^{j}:={\rm{support}}\left(\mu_{\sigma}^{j}\right)\subseteq\mathbb{R}^{d}\quad\forall(j,\sigma)\in\llbracket J\rrbracket\times\llbracket s\rrbracket,
\label{defxisigmaj}    
\end{align} 
the Cartesian product 
\begin{align}
\bm{\mathcal{X}}:=\prod_{(j,\sigma)\in\llbracket J\rrbracket\times\llbracket s\rrbracket}\mathcal{X}_{\sigma}^{j}\subseteq\left(\mathbb{R}^{d}\right)^{\otimes Js},
    \label{defCartesianProductSpaceMSBP}
\end{align} 
and a ground cost $\bm{C}:\bm{\mathcal{X}}\mapsto\mathbb{R}_{\geq 0}$.

Let $\mathcal{M}\left(\boldsymbol{\mathcal{X}}\right)$ denote the manifold of probability measures on  $\boldsymbol{\mathcal{X}}$, and for a fixed pair $(j,\sigma)\in\llbracket J\rrbracket \times \llbracket s\rrbracket$, let
\begin{align}
\bm{\mathcal{X}}_{(-j,-\sigma)}&:=\prod_{(a,b)\in\llbracket J\rrbracket \times \llbracket s\rrbracket\setminus (j,\sigma)}\mathcal{X}^{a}_{b}. 
\label{DefxiX}   
\end{align}
For a fixed $\varepsilon > 0$, the multimarginal Schr\"{o}dinger bridge problem (MSBP), generalizes the bimarginal problem \eqref{SBP} as
\begin{subequations}
\begin{align}
&\underset{\bm{M}\in \mathcal{M}\left(\bm{\mathcal{X}}\right)}{\min} \int_{\bm{\mathcal{X}}}\big\{\bm{C}\left(\bm{\xi}^{1}(\tau_{1}),\hdots,\bm{\xi}^{J}(\tau_{s})\right)\nonumber\\
&+ \varepsilon\log\bm{M}\left(\bm{\xi}^{1}(\tau_{1}),\hdots,\bm{\xi}^{J}(\tau_{s})\right)\}\differential\bm{M}\left(\bm{\xi}^{1}(\tau_{1}),\hdots,\bm{\xi}^{J}(\tau_{s})\right)\label{MSBPobj}\\
&\text{subject to}\int_{\bm{\mathcal{X}}_{(-j,-\sigma)}}\differential\bm{M}\left(\bm{\xi}^{1}(\tau_{1}), \hdots,\bm{\xi}^{J}(\tau_{s})\right) = \mu^{j}_{\sigma}\nonumber\\
&\qquad\qquad\qquad\qquad\qquad\qquad\forall(j,\sigma)\in\llbracket J\rrbracket\times\llbracket s\rrbracket.\label{MSBPconstr}
\end{align}
\label{MSBP}
\end{subequations}
The maximum likelihood interpretation for the bimarginal problem given in Proposition \ref{prop:SBPmaximumlikelihood} generalizes for MSBP \eqref{MSBP} as stated next in Proposition \ref{prop:MSBPmaximumlikelihood}. Similar to the case in Proposition \ref{prop:SBPmaximumlikelihood}, its proof is by direct computation and omitted.
\begin{proposition}\label{prop:MSBPmaximumlikelihood}
Given the sets \eqref{defxisigmaj}, \eqref{defCartesianProductSpaceMSBP}, ground cost $\bm{C}:\bm{\mathcal{X}}\mapsto\mathbb{R}_{\geq 0}$ and fixed $\varepsilon>0$, define 
\begin{align}
\bm{K}\left(\bm{\xi}^{1}(\tau_{1}),\hdots,\bm{\xi}^{J}(\tau_{s})\right) := \exp\left(-\dfrac{\bm{C}\left(\bm{\xi}^{1}(\tau_{1}),\hdots,\bm{\xi}^{J}(\tau_{s})\right)}{\varepsilon}\right),
\label{defkMulti}
\end{align} 
and the multimarginal Gibbs kernel 
\begin{align}
\boldsymbol{K}\left(\boldsymbol{\xi}^{1}\left(\tau_1\right), \ldots, \boldsymbol{\xi}^{J}\left(\tau_s\right)\right) \mu_1^{1} \otimes \ldots \otimes \mu^{J}_s.
\label{MultimarginalGibbsKernel}    
\end{align}
Let $\Pi\left(\mu^{1}_1,\hdots,\mu_{s}^{J}\right)$ denote the collection of measure-valued paths on $\mathcal{C}\left([\tau_1,\tau_s],\mathbb{R}^{d}\right)$ with $(j,\sigma)$ marginal $\mu_{\sigma}^{j}$ $\forall(j,\sigma)\in\llbracket J\rrbracket \times \llbracket s\rrbracket$. Then the optimal coupling $\bm{M}_{\rm{opt}}\in\mathcal{M}\left(\bm{\mathcal{X}}\right)$ in \eqref{MSBP} corresponds to the optimal measure-valued path $\pi_{\rm{opt}}\in\Pi\left(\mu^{1}_1,\hdots,\mu_{s}^{J}\right)$ solving the (scaled) relative entropy minimization problem:
\begin{align}
\underset{\pi\in\Pi(\mu^{1}_1,\hdots,\mu_{s}^{J})}{\min}\varepsilon D_{\rm{KL}}\left(\pi\Vert \bm{K}\left(\bm{\xi}^{1}(\tau_1),\hdots,\bm{\xi}^{J}(\tau_s)\right)\mu^{1}_1\otimes\hdots\otimes\mu^{J}_{s}\right).
\label{multimarginalKL}
\end{align}    
\end{proposition}
The existence-uniqueness of the minimizer for \eqref{multimarginalKL}, and thus for \eqref{MSBP}, follows from the same strict convexity argument as in Remark \ref{Remark:ExistenceUniquenessBimarginal}. 

Motivated by the maximum likelihood interpretation \eqref{multimarginalKL} for \eqref{MSBP}, we propose to solve \eqref{MSBP} for learning the stochastic computational resource state \eqref{ProblemStatementInformal}. The minimizer of \eqref{MSBP}, $\bm{M}^{\rm{opt}}\left(\bm{\xi}^{1}(\tau_{1}), \hdots,\bm{\xi}^{J}(\tau_{s})\right)$ can be used to compute the \emph{optimal} bimarginal probability mass transport plans between any $(j_1,\sigma_1),(j_2,\sigma_2)\in\llbracket J\rrbracket\times\llbracket s\rrbracket$, expressed as the bimarginal analogue of \eqref{MSBPconstr}:
\begin{equation}
    \int_{\bm{\mathcal{X}}_{(-j_1,-j_2,-\sigma_1,-\sigma_2)}}\differential\bm{M}^{\rm{opt}}\left(\bm{\xi}^{1}(\tau_{1}), \hdots,\bm{\xi}^{J}(\tau_{s})\right)
    \label{OptBimarginalCoupling}
\end{equation}
where
$$ \bm{\mathcal{X}}_{(-j_1,-j_2,-\sigma_1,-\sigma_2)}:=\prod_{(a,b)\in\llbracket J\rrbracket \times \llbracket s\rrbracket\setminus \{(j_1,\sigma_1)\cup(j_2,\sigma_2)\}}\mathcal{X}^{a}_{b}.$$
The coupling \eqref{OptBimarginalCoupling} will find use in implementing our multimarginal Sinkhorn recursions in Sec. \ref{sec:Algorithms}.

We next discuss the discrete formulations of \eqref{MSBP} that arise from the information graph structures in single and multi-core computing hardware for \emph{finite scattered data} $\{\bm{\xi}^{i,j}(\tau_{\sigma})\}_{i=1}^{n}$ and $\{\mu_{\sigma}^{j}\}_{(j,\sigma)\in\llbracket J\rrbracket\times\llbracket s\rrbracket}$ as in \eqref{DefEmpiricalMeasure}. We will see that the corresponding formulations lead to strictly convex problems over finite dimensional nonnegative tensors. To reduce the notational overload, we use the same boldfaced symbols for the continuous and discrete version of the tensors.

\subsection{Path Structured MSBP}\label{subsec:PathStructuredMSBP}
\begin{figure}[t]
    \centering    \includegraphics[width=0.99\linewidth]{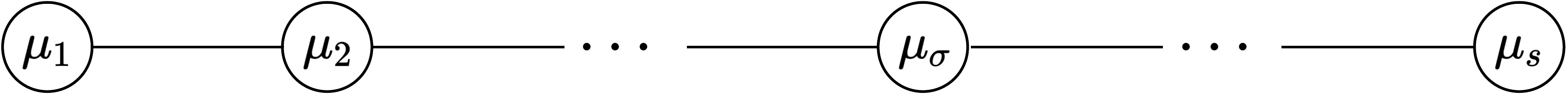}
    \caption{{\small{The path tree for sequentially observed $\{\mu_{\sigma}\}_{\sigma\in\llbracket s\rrbracket}$.}}}
\vspace*{-0.15in}
\label{FigPathTree}
\end{figure}
For single core ($J=1$) computing hardware, the stochastic process $\bm{\xi}(\tau)$ is indexed only by time $\tau\in[0,t]$. So the information graph structure is also induced by time, i.e., by sequential in time measure-valued observations $\{\mu_{\sigma}\}_{\sigma\in\llbracket s\rrbracket}$. In other words, the information graph in this case is a \emph{path tree} shown in Fig. \ref{FigPathTree}. We refer to an MSBP specialized to such path tree as the \emph{path structured MSBP}. 

Thanks to the path tree structure, the ground cost $\bm{C}$ in this case can be written as
\begin{align}
\bm{C}\!\left(\bm{\xi}(\tau_{1}),\hdots,\bm{\xi}(\tau_{s})\right) = \sum_{\sigma=1}^{s-1}c_{\sigma}\left(\bm{\xi}(\tau_{\sigma}),\bm{\xi}(\tau_{\sigma+1})\right)
\label{GroundCostPath}
\end{align}
where we choose the squared Euclidean distance sequential cost between two consecutive snapshot indices, i.e., $c_{\sigma}(\cdot,\cdot):=\|\cdot - \cdot\|_2^2$ $\forall\sigma\in\llbracket s\rrbracket$.

To formulate the \emph{discrete version} of the corresponding MSBP, notice that the ground cost in \eqref{MSBPobj}, in general, is an order $s$ tensor $\bm{C}\in\left(\mathbb{R}^{n}\right)^{\otimes s}_{\geq 0}$ having components
\begin{align}
\left[\bm{C}_{i_{1},\hdots,i_{s}}\right] = \bm{C}\left(\bm{\xi}_{i_1},\hdots,\bm{\xi}_{i_{s}}\right)
\label{CostTupleGeneric}
\end{align}
that encodes the cost of transporting unit amount of mass for an $s$ tuple $(i_1,\hdots,i_s)$. However, the path structured cost \eqref{GroundCostPath} implies that the $s$ tuple in \eqref{CostTupleGeneric} equals to $\sum_{\sigma=1}^{s-1}c_{\sigma}\left(\bm{\xi}_{i_{\sigma}},\bm{\xi}_{i_{\sigma+1}}\right)$, which is a sum of suitable elements of $s-1$ different (in our case, Euclidean) distance matrices.

The discrete mass tensor $\bm{M}\in\left(\mathbb{R}^{n}\right)^{\otimes s}_{\geq 0}$ has components 
\begin{align}
\left[\bm{M}_{i_{1},\hdots,i_{s}}\right] = \bm{M}\left(\bm{\xi}_{i_1},\hdots,\bm{\xi}_{i_{s}}\right),
\label{MassTupleGeneric}
\end{align}
where $\left[\bm{M}_{i_{1},\hdots,i_{s}}\right]$ denotes the amount of transported mass for an $s$ tuple $(i_1,\hdots,i_s)$.

Furthermore, for path structured MSBP, we can drop the core index superscript $j$ to simplify constraints \eqref{MSBPconstr} as
\begin{align}
\int_{\bm{\mathcal{X}}_{-\sigma}}\differential\bm{M}\left(\bm{\xi}(\tau_{1}), \hdots,\bm{\xi}(\tau_{s})\right) = \mu_{\sigma} \quad \forall\sigma\in\llbracket s\rrbracket.
\label{MSBPconstrPath}
\end{align}
In the discrete version, the empirical marginals $\bm{\mu}_{\sigma}\in\mathbb{R}^{n}_{\geq 0}$ are supported on the finite sets $\{\bm{\xi}^{i}(\tau_{\sigma})\}_{i=1}^{n}$ $\forall \sigma\in\llbracket s\rrbracket$. Then, the LHS of the constraints \eqref{MSBPconstrPath} are projections of the mass tensor \eqref{MassTupleGeneric}  on the $\sigma$th marginal $\bm{\mu}_{\sigma}$ $\forall \sigma\in\llbracket s\rrbracket$. We denote this projection as $\proj_{\sigma}(\bm{M})$, which is a mapping $\proj_{\sigma}:\left(\mathbb{R}^{n}\right)^{\otimes s}_{\geq 0}\mapsto \mathbb{R}^{n}_{\geq 0}$, and is given componentwise as
\begin{align}
\Big[{\rm{proj}}_\sigma \!\left(\bm{M}\right)_{r}\!\Big] \!\!=\!\! \!\!\sum_{i_1,\hdots,i_{\sigma-1},i_{\sigma + 1},\hdots, i_s}\!\!\!\!\!\!\!\!\bm{M}_{i_1,\hdots,i_{\sigma-1},r,i_{\sigma + 1},\hdots, i_s}.
\label{DefprojsigmaComponent}
\end{align} 

Similarly, the discrete version of \eqref{OptBimarginalCoupling} in the path structured case, is the projection of $\bm{M}\in\left(\mathbb{R}^{n}\right)^{\otimes s}_{\geq 0}$ on the $(\sigma_1,\sigma_2)$th marginal denoted as $\proj_{\sigma_1,\sigma_2}(\bm{M})$, i.e., $\proj_{\sigma_1,\sigma_2}:\left(\mathbb{R}^{n}\right)^{\otimes s}_{\geq 0}\mapsto \mathbb{R}^{n\times n}_{\geq 0}$, and is given componentwise as
\begin{align}
&\left[{\rm{proj}}_{\sigma_{1},\sigma_{2}} \!\left(\bm{M}\right)_{r,\ell}\!\right]\nonumber\\
&=\!\!\!\!\sum_{i_{\sigma}\mid\sigma\in\llbracket s\rrbracket\setminus\{\sigma_1,\sigma_2\}}\!\!\!\!\!\!\!\!\bm{M}_{i_1,\hdots,i_{\sigma_{1}-1},r,i_{\sigma_{1} + 1},\hdots,i_{\sigma_{2}-1},\ell,i_{\sigma_{2} + 1},\hdots,i_s}.
\label{DefprojsigmaComponentDouble}
\end{align}
Then, the discrete path structured version of \eqref{MSBP} becomes
\begin{subequations}
\begin{align}
&\underset{\bm{M}\in\left(\mathbb{R}^{n}\right)^{\otimes s}_{\geq 0}}{\min}~\langle\bm{C}+\varepsilon\log\bm{M},\bm{M}\rangle\label{DiscreteMSBPobj}\\
&\text{subject to}~~{\rm{proj}}_{\sigma}\left(\bm{M}\right) = \bm{\mu}_{\sigma}\quad\forall\sigma\in\llbracket s\rrbracket.\label{DiscereteMSBPconstr}
\end{align}
\label{DiscreteMSBP}
\end{subequations}
Notice that \eqref{DiscreteMSBP} is a strictly convex problem in $n^{s}$ decision variables, and is computationally intractable in general. In Sec. \ref{sec:Algorithms}, we will design algorithms to solve them efficiently.


\subsection{Barycentric MSBP}\label{subsec:BaryStructuredMSBP}
\begin{figure}[t]
    \centering    \includegraphics[width=0.99\linewidth]{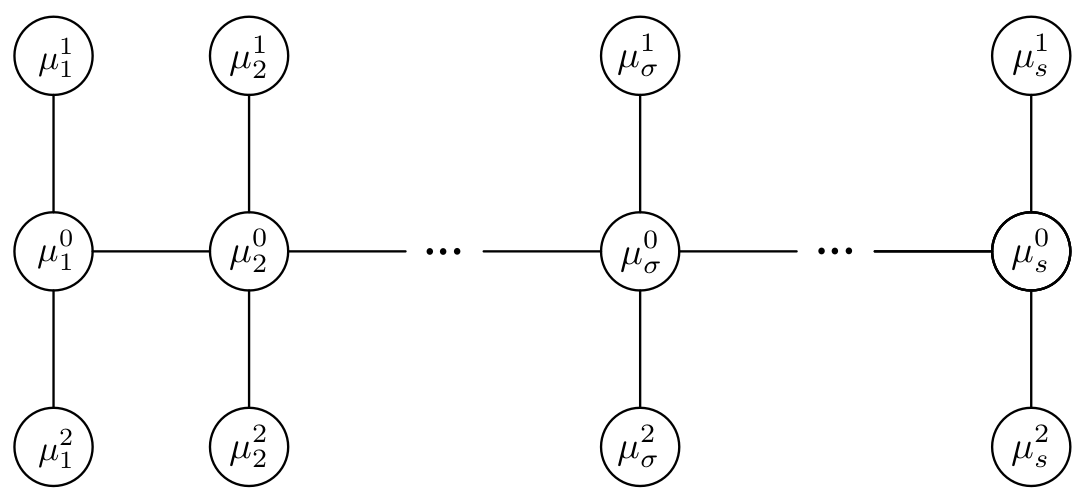}
    \caption{{\small{The tree graph for barycentric MSBP (Sec. \ref{subsec:BaryStructuredMSBP}) with measures $\{\mu_{\sigma}^j\}_{(j,\sigma)\in\left(\{0\}\cup\llbracket J\rrbracket\right)\times\llbracket s\rrbracket}$, where $J=2$.}}}
\vspace*{-0.15in}
\label{FigBC_PathTree}
\end{figure}
For multi-core ($J>1$) computing hardware, the stochastic process $\bm{\xi}^{j}(\tau)$ is indexed by both time $\tau\in[0,t]$ and CPU core $j\in\llbracket J\rrbracket$. Thus, the information graph structure is induced jointly by time and cores, i.e., by measure-valued observations $\{\mu_{\sigma}^{j}\}_{(j,\sigma)\in\llbracket J\rrbracket\times\llbracket s\rrbracket}$. Unlike Sec. \ref{subsec:PathStructuredMSBP}, we now have $J$ measure-valued snapshots at each time index. To account for this, we propose two different MSBP formulations, the first of which termed \emph{barycentric MSBP} is discussed next.

The main idea behind our proposed barycentric MSBP formulation is to imagine a phantom CPU core whose statistics is the average (i.e., barycenter) of the measure-valued snapshots for all cores at any fixed time index $\sigma\in\llbracket s\rrbracket$. Accordingly, in this formulation, we consider the core index $j\in\{0\}\cup\llbracket J\rrbracket$ where $j=0$ refers to the phantom barycentric CPU. The corresponding graph structure for $J=2$ is shown in Fig. \ref{FigBC_PathTree}. While a barycentric graph as in Fig. \ref{FigBC_PathTree} is more general than the path tree as in Fig. \ref{FigPathTree}, we note that the \emph{treewidth} \cite[p. 354-355]{diestel2005graph} for barycentric graphs equals unity as in the path tree case. It is known that the computational complexity for graph-structured MSBP problems grow with the treewidth \cite{fan2022complexity,altschuler2023polynomial}. We will discuss the specific complexity for the proposed algorithms in Sec. \ref{sec:Algorithms}.

The barycentric graph structure implies that \eqref{CostTupleGeneric} is no longer equal to $\sum_{\sigma=1}^{s-1}c_{\sigma}\left(\bm{\xi}_{i_{\sigma}},\bm{\xi}_{i_{\sigma+1}}\right)$ as in Sec. \ref{subsec:PathStructuredMSBP}, but instead becomes a sum of two types of ground transport costs, viz. the cost of transport between consecutive-in-time barycenters and the cost of transport between barycentric and actual CPU cores. Denoting the barycentric random vectors as $\{\bm{\xi}^{0}(\tau_{\sigma})\}_{\sigma\in\llbracket s\rrbracket}$, and letting $c_{j,\sigma}(\cdot,\cdot)$ be the corresponding ground transport costs for all $j\in\{0\}\times\llbracket J\rrbracket$, the cost \eqref{CostTupleGeneric} for barycentric MSBP equals
$$\sum_{\sigma=1}^{s-1}c_{0,\sigma}\left(\bm{\xi}_{i_{\sigma}}^{0},\bm{\xi}^{0}_{i_{\sigma+1}}\right) + \sum_{\sigma=1}^{s}\sum_{j=1}^{J}c_{j,\sigma}\left(\bm{\xi}^{j}_{i_{\sigma}},\bm{\xi}^{0}_{i_{\sigma}}\right).$$
Similar ideas appeared in \cite[Sec. 3.3]{elvander2020multi} in a different context.

Defining the barycentric index set 
\begin{align}
\Lambda_{\rm{BC}}:=(\{0\}\cup\llbracket J\rrbracket)\times\llbracket s\rrbracket,
    \label{BaryIndexSet}
\end{align}
we note that $\vert\Lambda_{\rm{BC}}\vert = (J+1)s$, and thus for the barycentric MSBP, the tensors $\bm{C},\bm{M}\in\left(\mathbb{R}^{n}\right)^{\otimes (J+1)s}_{\geq 0}$. In summary, the proposed barycentric MSBP is
\begin{subequations}
\begin{align}
&\underset{\bm{M}\in\left(\mathbb{R}^{n}\right)^{\otimes (J+1)s}_{\geq 0}}{\min}~\langle\bm{C}+\varepsilon\log\bm{M},\bm{M}\rangle\label{DiscreteBaryMSBPobj}\\
&\text{subject to}~~{\rm{proj}}_{(j,\sigma)}\left(\bm{M}\right) = \bm{\mu}_{\sigma}^{j}\quad\forall(j,\sigma)\in\Lambda_{\rm{BC}}.\label{DiscereteBaryMSBPconstr}
\end{align}
\label{DiscreteBaryMSBP}
\end{subequations}
In Sec. \ref{subsec:BaryProjections}, we will discuss the computation of projections \eqref{MSBPconstr} and \eqref{OptBimarginalCoupling} for the minimizer of \eqref{DiscreteBaryMSBP}.
 

\subsection{Series-Parallel Graph Structured MSBP}\label{subsec:SeriesParallelGraphStructuredMSBP}
\begin{figure}[t]
    \centering    \includegraphics[width=0.99\linewidth]{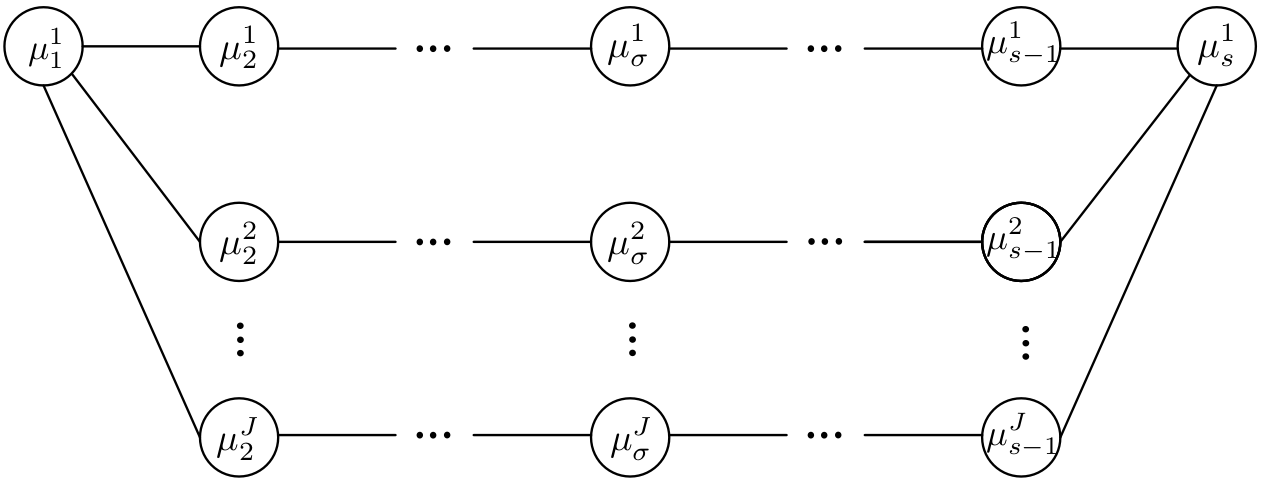}
    \caption{{\small{The graph for series-parallel graph-structured MSBP (Sec. \ref{subsec:SeriesParallelGraphStructuredMSBP}) with $J$ parallel paths of length $s$.}}}
\vspace*{-0.15in}
\label{FigSP_PathTree}
\end{figure}
Different from the barycentric MSBP in Sec. \ref{subsec:BaryStructuredMSBP}, we now propose another MSBP formulation for the multi-core ($J>1$) case, based on the series-parallel information graph structure in Fig. \ref{FigSP_PathTree}. 

This \emph{series-parallel graph structured MSBP} is motivated by the observation that many software on multi-core computing hardware have notions of input and output \emph{terminals}.  In such applications, the input and output of the computational task are aggregated to a single core. Thus, the information graph structure is induced by measure-valued observations $\{\mu_{\sigma}^{j}\}_{(j,\sigma)}$ where the tuple $(j,\sigma)$ belongs to the index set 
\begin{align}
\Lambda_{\rm{SP}} := \left(\llbracket J\rrbracket\times\llbracket s\rrbracket\right) \setminus \left(\left(\llbracket J\rrbracket\setminus\{1\}\right)\times\{1,s\}\right).
\label{SPIndexSet}
\end{align}
Although the series-parallel graph in Fig. \ref{FigSP_PathTree} is not a tree, it has treewidth at most two. So such MSBPs remain computationally tractable as before (see details in Sec. \ref{sec:Algorithms}).

Unlike the formulation in Sec. \ref{subsec:BaryStructuredMSBP}, the series-parallel graph structured MSBP formulation does not involve any phantom CPU core. From Fig. \ref{FigSP_PathTree}, the cost \eqref{CostTupleGeneric} now equals
\begin{align*}
&\sum_{j=1}^{J}\bigg\{c_{j,1}\left(\bm{\xi}_{i_{1}}^{j},\bm{\xi}_{i_{2}}^{j}\right) + c_{j,s-1}\left(\bm{\xi}_{i_{s-1}}^{j},\bm{\xi}_{i_{s}}^{j}\right)\bigg\}\\
&\qquad+\sum_{\sigma=2}^{s-1}\sum_{j=1}^{J}c_{j,\sigma}\left(\bm{\xi}_{i_{\sigma}}^{j},\bm{\xi}_{i_{\sigma+1}}^{j}\right).
\end{align*}

By Lemma \ref{LemmaCardnilaitySetMinus}, we find $\vert \Lambda_{\rm{SP}}\vert = J(s-2) + 2$, and hence for the series-parallel graph structured MSBP, the tensors $\bm{C},\bm{M}\in\left(\mathbb{R}^{n}\right)^{\otimes J(s-2) + 2}_{\geq 0}$. In summary, the proposed series-parallel graph-structured MSBP is
\begin{subequations}
\begin{align}
&\underset{\bm{M}\in\left(\mathbb{R}^{n}\right)^{\otimes (J(s-2)+2)}_{\geq 0}}{\min}~\langle\bm{C}+\varepsilon\log\bm{M},\bm{M}\rangle\label{DiscreteSPMSBPobj}\\
&\text{subject to}~~{\rm{proj}}_{(j,\sigma)}\left(\bm{M}\right) = \bm{\mu}_{\sigma}^{j}\quad\forall(j,\sigma)\in\Lambda_{\rm{SP}}.\label{DiscereteSPMSBPconstr}
\end{align}
\label{DiscreteSPMSBP}
\end{subequations}
In Sec. \ref{subsec:SPProjections}, we will discuss the computation of projections \eqref{MSBPconstr} and \eqref{OptBimarginalCoupling} for the minimizer of \eqref{DiscreteSPMSBP}.

\begin{table}[t]
\centering
\begin{tabular}{ | c | c | c | } 
\hline
Graph structure & MSBP & Index set $\Lambda$ \\
\hline\hline
Path tree (Sec. \ref{subsec:PathStructuredMSBP}) & \eqref{DiscreteMSBP} & $\llbracket s\rrbracket$\\
& & \\
Barycentric (Sec. \ref{subsec:BaryStructuredMSBP}) & \eqref{DiscreteBaryMSBP} & $\Lambda_{\rm{BC}}$ in \eqref{BaryIndexSet}\\
& &\\
Series-parallel (Sec. \ref{subsec:SeriesParallelGraphStructuredMSBP}) & \eqref{DiscreteSPMSBP} & $\Lambda_{\rm{SP}}$ in \eqref{SPIndexSet}\\ 
\hline
\end{tabular}
\caption{The index set $\Lambda$ in graph structured MSBPs.}
\label{TablefLambda}
\vspace*{-0.15in}
\end{table}
\section{Algorithms}\label{sec:Algorithms}
This section provides algorithmic details and computational complexities to solve the discrete MSBPs \eqref{DiscreteMSBP}, \eqref{DiscreteBaryMSBP} and \eqref{DiscreteSPMSBP}. 

These MSBPs are strictly convex tensor optimization problems in $n^{s}, n^{(J+1)s},$ and $n^{J(s-2)+2}$ decision variables, respectively, and computationally intractable in general. By leveraging an interplay between duality and graph structures, we will see that it is possible to reduce the computational complexity from exponential to linear in $s$.


Recognizing that \eqref{DiscreteMSBP}, \eqref{DiscreteBaryMSBP}, \eqref{DiscreteSPMSBP} are instances of the generic structured optimization problem:
\begin{subequations}
\begin{align}
&\underset{\bm{M}\in\left(\mathbb{R}^{n}\right)^{\otimes \vert \Lambda\vert}_{\geq 0}}{\min}~\langle\bm{C}+\varepsilon\log\bm{M},\bm{M}\rangle\label{GenericDiscreteMSBPobj}\\
&\text{subject to}~~{\rm{proj}}_{(j,\sigma)}\left(\bm{M}\right) = \bm{\mu}_{\sigma}^{j}\quad\forall(j,\sigma)\in\Lambda,\label{GenericDiscereteMSBPconstr}
\end{align}
\label{GenericDiscreteMSBP}
\end{subequations}
for suitable index set $\Lambda$ (see Table \ref{TablefLambda}), the following is consequence of strong Lagrange duality.
\begin{proposition}\cite{benamou2015iterative,marino2020optimal,carlier2022linear}
\label{propSinkhornConvergence}
Given problem \eqref{GenericDiscreteMSBP}, let $\bm{\lambda}^{j}_{\sigma}\in\mathbb{R}^{n}$ be the Lagrange multipliers for the equality constraints \eqref{GenericDiscereteMSBPconstr} for all $(j,\sigma)\in\Lambda$. Let 
\begin{subequations}
\begin{align}
    \bm{K} &:= \exp(-\bm{C}/\varepsilon)\in\left(\mathbb{R}^{n}\right)^{\otimes \vert \Lambda\vert}_{> 0}, \label{eqstructKtensor} \\
    \bm{u}^{j}_{\sigma} &:= \exp(\bm{\lambda}^{j}_{\sigma}/\varepsilon)\in\mathbb{R}^{n}_{>0}\quad\forall(j,\sigma)\in\Lambda, \label{eqstructuvector}\\
    \bm{U} &:= \displaystyle\otimes_{(j,\sigma)\in\Lambda}\bm{u}^{j}_{\sigma}\in\left(\mathbb{R}^{n}\right)^{\otimes \vert \Lambda\vert}_{> 0}. \label{eqstructUtensor}
\end{align}
\label{eqstructKUutensor}
\end{subequations}
Then, the multi-marginal Sinkhorn recursions
\begin{align}
\bm{u}^{j}_{\sigma} \leftarrow \bm{u}^{j}_{\sigma} \odot \bm{\mu}^{j}_{\sigma}\oslash{\rm{proj}}_{(j,\sigma)}\left(\bm{K}\odot\bm{U}\right) ~\forall(j,\sigma)\in\Lambda,
\label{MultimarginalSink}    
\end{align}
have guaranteed linear rate of convergence\footnote{The proof sketch for discrete state space is as follows. The MSBP \eqref{GenericDiscreteMSBP} can be recast \cite[Sec. 4.1]{benamou2015iterative} as a Kullback-Leibler projection to a convex set that is an intersection of $|\Lambda|$ hyperplanes given by \eqref{GenericDiscereteMSBPconstr}. That this iterative Kullback-Leibler projection has guaranteed convergence follows from the seminal result on iterative Bregman projection \cite{bregman1967relaxation}, and from the fact that the Kullback-Leibler divergence is an instance of Bregman divergence. That the rate is linear follows from the co-ordinate descent analysis by Luo and Tseng \cite[p. 25-26]{luo1992convergence}.

For the continuous state space, ref. \cite[Thm. 4.7]{marino2020optimal} proves the convergence of multi-marginal recursions. Ref. \cite[Sec. 3]{carlier2022linear} extends the result of \cite{marino2020optimal} by showing linear rate of convergence.}, and the minimizer $\bm{M}^{\rm{opt}}$ for \eqref{GenericDiscreteMSBP} is given in terms of the converged $\bm{U}$ as
\begin{align}
\bm{M}^{\rm{opt}}=\bm{K}\odot\bm{U}.
    \label{Mopt}
\end{align}
\end{proposition}
Once $\bm{M}^{\rm{opt}}$ is computed, its bimarginal projections of the form \eqref{DefprojsigmaComponentDouble} yield the probability mass transport matrices between any two marginals. With this, we are able to interpolate between any two marginals to obtain a predicted distribution. 

\begin{remark}\label{RemarkInterpolation}
    In our problem, the useful application of this interpolation is as follows: given a time $\tau\in[0,t)$ and CPU index $j\in\llbracket J\rrbracket$, find $\sigma\in\llbracket s\rrbracket$ such that $\tau_\sigma\leq\tau<\tau_{\sigma+1}$, and use the expressions for the bimarginal projections above to compute the bimarginal bridge
    \begin{equation*}
        M^{j,\sigma}:={\rm{proj}}_{(j,\sigma),(j,\sigma+1)}(\bm{M}^{\rm{opt}}) : \mu_\sigma^j \to \mu_{\sigma+1}^j \quad\left(\in\R^{n\times n}_{\geq 0}\right).
    \end{equation*}
    Using this projection, we can interpolate between $\mu_\sigma^j$ and $\mu_{\sigma+1}^j$ to obtain our estimate $\widehat{\mu}_{\tau}^j$ for the computational resource usage distribution for CPU $j$ at time $\tau$, as
    \begin{align}
    \!\hat{\mu}_\tau^j \!:=\!\sum_{r=1}^n\!\sum_{\ell=1}^n \!\left[\!M^{j,\sigma}_{r,\ell}\!\right]\!\!\delta(\bm{\xi}^j-\widehat{\bm{\xi}}^j(\tau,\bm{\xi}^{r,j}(\tau_{\sigma}),\bm{\xi}^{\ell,j}(\tau_{\sigma+1})))
    \label{MuInterpolation}
    \end{align}
    where $\widehat{\bm{\xi}}^j(\tau,\bm{\xi}^{r,j}(\tau_{\sigma}),\bm{\xi}^{\ell,j}(\tau_{\sigma+1}))\!\!:= \!(1-\lambda)\bm{\xi}^{r,j}(\tau_{\sigma}) \!+\! \lambda\bm{\xi}^{\ell,j}(\tau_{\sigma+1})$,
    and $\lambda:=\dfrac{\tau-\tau_{\sigma}}{\tau_{\sigma+1}-\tau_{\sigma}}\in[0,1]$.
\end{remark}

For numerically solving bi-marginal SBPs (i.e., the case $s=2, J=1$), the Sinknorn recursions \eqref{MultimarginalSink} have become the standard \cite{cuturi2013sinkhorn,genevay2018learning,flamary2021pot,10347388}. However, applying the same is challenging for MSBPs because computing $\proj_{(j,\sigma)}\left(\bm{K}\odot\bm{U}\right)$ requires $\mathcal{O}\left(n^{\vert\Lambda\vert}\right)$ operations. The same issue arises for computing $\proj_{(j_1,\sigma_1),(j_2,\sigma_2)}\left(\bm{K}\odot\bm{U}\right)$.
For both projections, this complexity can be reduced by exploiting the structure of the Hilbert-Schmidt inner product $\langle\bm{K},\bm{U}\rangle$. Specifically, we make use of the following results from \cite{elvander2020multi}.

\begin{lemma}\label{lemmaproj1and2}
Let the tensor $\bm{U}$ be as in \eqref{eqstructUtensor}, and consider a tensor $\bm{K}\in\left(\mathbb{R}^{n}\right)^{\otimes \vert\Lambda\vert}$.\\
\noindent (i) \cite[Lemma 1]{elvander2020multi} For a fixed $(j,\sigma)\in\Lambda$, if $\langle\bm{K},\bm{U}\rangle = \bm{w}_1^\top\diag\left(\bm{u}^{j}_{\sigma}\right)\bm{w}_2$ for some vectors $\bm{w}_1,\bm{w}_2\in\R^n$ that do not depend on $\bm{u}^{j}_{\sigma}$, then
$${\rm{proj}}_{(j,\sigma)}(\bm{K}\odot\bm{U})=\bm{w}_1\odot\bm{u}^{j}_\sigma\odot\bm{w}_2.$$
\noindent (ii) \cite[Lemma 2]{elvander2020multi} For fixed $(j_1,\sigma_1),(j_2,\sigma_2)\in\Lambda$, if $\langle\bm{K},\bm{U}\rangle = \bm{w}_1^\top\diag(\bm{u}^{j_{1}}_{\sigma_1})\Phi\diag(\bm{u}^{j_{2}}_{\sigma_2})\bm{w}_3$ for some vectors $\bm{w}_1,\bm{w}_3\in\R^n$ and matrix $\Phi\in\R^{n\times n}$ where $\bm{w}_1,\Phi,\bm{w}_{3}$ do not depend on $\bm{u}^{j_{1}}_{\sigma_1},\bm{u}^{j_{2}}_{\sigma_2}$, then
$${\rm{proj}}_{(j_1,\sigma_1),(j_2,\sigma_2)}(\bm{K}\odot\bm{U}) = \diag(\bm{w}_1\odot\bm{u}^{j_1}_{\sigma_1})\Phi\diag(\bm{u}^{j_2}_{\sigma_2}\odot\bm{w}_3).$$
\end{lemma}

The special inner product structures in Lemma \ref{lemmaproj1and2} arises from structured tensors $\bm{K}$ which are in turn induced by the graph structures discussed earlier. Thus, expressing $\langle\bm{K},\bm{U}\rangle$ in the appropriate forms helps compute the desired projections. Below, we show how the imposition of structure on $\bm{K}$ allows this to be done efficiently.



\subsection{Projections for Path Structured MSBP}\label{subsec:ProjectionForPath}
For software running on a single CPU core, we have $J=1$ as in Sec. \ref{subsec:PathStructuredMSBP}. The corresponding MSBP is \eqref{DiscreteMSBP}. Here, the cost tensor $\bm{C}$ in \eqref{CostTupleGeneric} has a path structure 
\begin{equation}
    \left[\bm{C}_{(i_\sigma \mid \sigma\in\llbracket s \rrbracket)}\right] = \left[\bm{C}_{i_1,\dots,i_s}\right] = \sum_{\sigma=1}^{s-1}\left[C_{i_{\sigma},i_{\sigma+1}}^{\sigma}\right]
\label{eqC_SingleCore}
\end{equation}
where the matrix $C^{\sigma}\in\mathbb{R}^{n\times n}_{\geq 0}$ encodes the cost of transporting unit mass from $\{\bm{\xi}^{i}(\tau_{\sigma})\}_{i=1}^{n}$ to $\{\bm{\xi}^{i}(\tau_{\sigma+1})\}_{i=1}^{n}$. This allows us to write $\bm{K}$ in \eqref{eqstructKtensor} as
\begin{equation}
    \left[\bm{K}_{(i_\sigma \mid \sigma\in\llbracket s \rrbracket)}\right] = \prod_{\sigma=1}^{s-1}\left[K_{i_{\sigma},i_{\sigma+1}}^{\sigma}\right] 
\label{eqK_SingleCore}
\end{equation}
which leads to the following expressions for the marginal projections.

\begin{proposition}\label{PropMultiSinkPathCost} \cite[Prop. 2]{elvander2020multi}, \cite[Prop. 1]{bondar2023path}
If $\bm{C}$ has the form \eqref{eqC_SingleCore}, $K^{\sigma}:=\exp(-C^{\sigma}/\varepsilon)\in\mathbb{R}^{n\times n}_{\geq 0}$, $\bm{K}$ as in \eqref{eqK_SingleCore}, and $\bm{U}$ as in \eqref{eqstructUtensor}, then \eqref{DefprojsigmaComponent} and \eqref{DefprojsigmaComponentDouble} can be expressed as
\begin{align}
&{\rm{proj}}_{\sigma}\!\left(\bm{K}\odot\bm{U}\right)\!=\!\!\left(\!\bm{u}_1^{\!\top} K^{1}\!\prod_{k=2}^{\sigma-1}\!\diag(\bm{u}_k)K^{k}\!\!\right)^{\!\!\top}\!\!\!\!\odot\bm{u}_{\sigma}\odot \nonumber\\
&\left(\!\!\left(\prod_{k=\sigma+1}^{s-1}K^{k-1}\diag(\bm{u}_k)\!\right)\!\!K^{s-1}\bm{u}_s\!\right)\,\forall\sigma\in\llbracket s\rrbracket,\label{ProjSimplified}
\end{align}
and ${\rm{proj}}_{\sigma_1,\sigma_2}\!\left(\bm{K}\odot\bm{U}\right)=$
\begin{align}
&\diag\!\left(\!\!\bm{u}_1^\top K^{1}\prod_{k=2}^{\sigma_1-1}\diag(\bm{u}_k)K^{k}\!\!\right)\prod_{k=\sigma_1+1}^{\sigma_2}\!\!\left(K^{k-1}\diag(\bm{u}_k)\right)\nonumber\\
&\diag\!\left(\!\!\left(\prod_{k=\sigma_2+1}^{s-1}K^{k-1}\diag(\bm{u}_k)\!\!\right)\!K^{s-1}\bm{u}_s\!\!\right)\nonumber\\
&\qquad\qquad\qquad\qquad~\forall(\sigma_1,\sigma_2)\in\{\llbracket s\rrbracket^{\otimes 2}\mid\sigma_1 < \sigma_2\}.
\label{Proj2Simplified}
\end{align}
\end{proposition}

Observe that even the na\"ive computation of \eqref{ProjSimplified} is dominated by $2s-4$ matrix-vector multiplications (by cancellation, \eqref{MultimarginalSink} can be computed by $s-1$ such multiplications; see \cite[Remark 3]{bondar2023path}). Since such multiplications have $\mathcal{O}(n^2)$ complexity, each Sinkhorn iteration has $\mathcal{O}\left((s-1)n^2\right)$ complexity -- linear in $s$, and a great improvement from the general $\mathcal{O}(n^s)$ complexity of the method. Our method's $\mathcal{O}\left((s-1)n^2\right)$ complexity is sharp for exact computation. The recent work \cite{ba2022accelerating} reduces the complexity in $n$ from quadratic to linear at the expense of approximate computation. 

\begin{remark}\label{RemarkProj2Complexity}
    While it is clear from their expressions that the bimarginal projection \eqref{Proj2Simplified} has similar order-of-magnitude complexity to the unimarginal projection \eqref{ProjSimplified}, hereafter we focus only on the complexity of the latter, as these are performed every Sinkhorn iteration until the method converges. Following this, bimarginal projections of the solution $\bm{M}^{\rm{opt}}$ onto each pair of marginals of interest are to be performed only once a posteriori, and so exact floating-point operational count for these projections is not critical.
\end{remark}

\subsection{Projections for Barycentric MSBP}\label{subsec:BaryProjections}
For software running on multiple CPU cores, we have $J>1$. The corresponding barycentric MSBP \eqref{DiscreteBaryMSBP} formulated in Sec. \ref{subsec:BaryStructuredMSBP} involves a tree that is more general than a path. 

In this subsection, we let $n_{0}$ denote the number of samples in the barycentric CPU. For the sake of generality, here we derive the complexity for the projections in terms of $n$ and $n_{0}$. Then, $\bm{C},\bm{M},\bm{U}\in\left(\mathbb{R}^{n}\right)^{\otimes Js}_{\geq 0}\otimes \left(\mathbb{R}^{n_{0}}\right)^{\otimes s}_{\geq 0}$. For $n_0 \neq n$, formulation \eqref{GenericDiscreteMSBP}-\eqref{eqstructKUutensor} applies mutatis-mutandis. 

Recall that here the index set $\Lambda=\Lambda_{\rm{BC}}$ as in \eqref{BaryIndexSet}. The cost tensor $\bm{C}$ for barycentric MSBP takes the form 
\begin{align}
    \left[\bm{C}_{(i_{(j,\sigma)} \mid (j,\sigma)\in\Lambda_{\rm{BC}})}\right] = \sum_{\sigma=1}^{s-1}\left[C_{i_{(0,\sigma)},i_{(0,\sigma+1)}}^{0,\sigma}\right] \nonumber \\
    + \sum_{\sigma=1}^{s}\sum_{j=1}^{J}\left[C_{i_{(j,\sigma)},i_{(0,\sigma)}}^{j,\sigma}\right]
\label{eqC_BC}
\end{align}
where the matrices $C^{0,\sigma}\in\R^{n_{0}\times n_{0}}_{\geq 0}$ are the ground cost matrices between barycenters $\bm{\xi}^0(\tau_\sigma)$ and $\bm{\xi}^0(\tau_{\sigma+1})$ for $\sigma\in\llbracket s-1\rrbracket$, whereas $C^{j,\sigma}\in\R^{n_{0}\times{n}}_{\geq 0}$ are those between the CPU marginals $\bm{\xi}^j(\tau_\sigma)$ and their barycenters $\bm{\xi}^0(\tau_\sigma)$, for $(j,\sigma)\in\llbracket J\rrbracket\times\llbracket s\rrbracket$. So now, we can write $\bm{K}$ in \eqref{eqstructKtensor} as
\begin{align}
&\left[\bm{K}_{(i_{(j,\sigma)} \mid (j,\sigma)\in\Lambda_{\rm{BC}})}\right] \nonumber\\
    &= \bigg(\prod_{\sigma=1}^{s-1}\left[K_{i_{(0,\sigma)},i_{(0,\sigma+1)}}^{0,\sigma}\right]\bigg)
    \times \prod_{\sigma=1}^{s}\prod_{j=1}^{J}\left[K_{i_{(j,\sigma)},i_{(0,\sigma)}}^{j,\sigma}\right].
\label{eqK_BC}
\end{align}

With this, the projections can be computed as follows (proof in Appendix \ref{AppProofPropBCSinkPathCost}).

\begin{proposition}\label{PropBCSinkPathCost}
If $\bm{C}$ has the form \eqref{eqC_BC}, $K^{j,\sigma}:=\exp(-C^{j,\sigma}/\varepsilon)$, $\bm{K}$ as in \eqref{eqK_BC}, and $\bm{U}$ as in \eqref{eqstructUtensor}, then the projections \eqref{DefprojsigmaComponent} and \eqref{DefprojsigmaComponentDouble} can be expressed as
\begin{subequations}
\begin{align}
&{\rm{proj}}_{(0,\sigma)}\!\left(\bm{K}\odot\bm{U}\right)\! = \!\!\left(\!\bm{p}_1^{\!\top} K^{0,1}\!\prod_{k=2}^{\sigma-1}\!\diag(\bm{p}_k)K^{0,k}\!\!\right)^{\!\!\top}\!\!\!\!\odot\bm{p}_{\sigma}\odot \nonumber\\ &\quad\left(\!\!\left(\prod_{k=\sigma+1}^{s-1}K^{0,k-1}\diag(\bm{p}_k)\!\right)\!\!K^{0,s-1}\bm{p}_s\!\right) \quad\forall\sigma\in\llbracket s\rrbracket,\label{ProjSimplified_BC_marg0} \\   &{\rm{proj}}_{(j,\sigma)}\!\left(\bm{K}\odot\bm{U}\right)\! = \bm{u}_{\sigma}^{j}\odot {K^{j,\sigma}}^{\top}\!\Bigg(\!\!\left(\!\bm{p}_1^{\!\top} K^{0,1}\!\prod_{k=2}^{\sigma-1}\!\diag(\bm{p}_k)K^{0,k}\!\!\right)^{\!\!\top} \nonumber\\  &\odot\left(\bm{p}_{\sigma}\oslash\left({K^{j,\sigma}}\bm{u}_{\sigma}^{j}\right)\right)\odot\left(\!\!\left(\prod_{k=\sigma+1}^{s-1}K^{0,k-1}\diag(\bm{p}_k)\!\right)\!\!K^{0,s-1}\bm{p}_s\!\!\right)\nonumber\\
&\qquad\qquad\qquad\qquad\qquad\qquad\forall(j,\sigma)\in\llbracket J\rrbracket\times\llbracket s\rrbracket,\label{ProjSimplified_BC_marggen}
\end{align}
\label{ProjSimplified_BC}
\end{subequations}
and for $(j,\sigma)\in\llbracket J\rrbracket\times\llbracket s\rrbracket$ and $\sigma_1,\sigma_2\in\llbracket s\rrbracket$ such that $\sigma_1<\sigma_2$,
\begin{subequations}
\begin{align}
&{\rm{proj}}_{(0,\sigma),(j,\sigma)}\!\left(\bm{K}\odot\bm{U}\right) = \diag(\bm{u}_{\sigma}^{0})\diag\left({K^{0,\sigma-1}}^{\top}\bm{\rho}_{(0,\sigma),(j,\sigma)}\right) \nonumber \\
&\quad\quad\quad\quad\quad\quad\quad\quad\quad\quad\quad{K^{j,\sigma}}\diag(\bm{u}_{\sigma}^{j}),\label{Proj2Simplified_BC_marg0} \\
&{\rm{proj}}_{(0,\sigma_1),(0,\sigma_2)}\!\left(\bm{K}\odot\bm{U}\right) = \diag\!\left(\!\!\bm{p}_1^\top K^{0,1}\prod_{k=2}^{\sigma_1-1}\diag(\bm{p}_k)K^{0,k}\!\!\right)\nonumber\\
&\qquad\qquad\qquad\qquad~~~\odot\diag(\bm{p}_{\sigma_1})\!\!\prod_{k=\sigma_1+1}^{\sigma_2}\!\!\left(K^{0,k-1}\diag(\bm{p}_k)\right)\nonumber\\
&\qquad\qquad\qquad~\odot\diag\!\left(\!\!\left(\prod_{k=\sigma_2+1}^{s-1}K^{0,k-1}\diag(\bm{p}_k)\!\!\right)\!K^{0,s-1}\bm{p}_s\!\!\right)\label{Proj2Simplified_BC_margs}
\end{align}
\label{Proj2Simplified_BC}
\end{subequations}
where $\bm{p}_\sigma:=\bm{u}_{\sigma}^{0}\odot(\bigodot_{j\in\llbracket J\rrbracket}K^{j,\sigma}\bm{u}_{\sigma}^{j})$ for $\sigma\in\llbracket s\rrbracket$, and
\begin{align*}
    \bm{\rho}_{(0,\sigma),(j,\sigma)} &:= \!\!\left(\!\bm{p}_1^{\!\top} K^{0,1}\!\prod_{k=2}^{\sigma-1}\!\diag(\bm{p}_k)K^{0,k}\!\!\right)^{\!\!\top} \nonumber\\
    &\quad\odot\left( \bm{p}_{\sigma}\oslash\left(\bm{u}_{\sigma}^{0}\odot K^{j,\sigma}\bm{u}_{\sigma}^{j}\right) \right)\odot \nonumber\\
    &\left(\!\!\left(\prod_{k=\sigma+1}^{s-1}K^{0,k-1}\diag(\bm{p}_k)\!\right)\!\!K^{0,s-1}\bm{p}_s\!\right)\quad\forall\sigma\in\llbracket s\rrbracket.
\end{align*}
\end{proposition}
Similar to the single-core case in the Sec. \ref{subsec:ProjectionForPath}, the computational complexity for all projections in Proposition \ref{PropBCSinkPathCost} are linear in $J$ and $s$. Specifically, the computation of the $s$ $\bm{p}_{\sigma}$ vectors in Proposition \ref{PropBCSinkPathCost} requires a total of $Js$ matrix-vector multiplications, and the total number of floating-point operations when projecting onto a barycenter (as in \eqref{ProjSimplified_BC_marg0}) is
\begin{equation}
    Js(n_{0}n + n_{0}) + \left(2n_{0}\right) + (2s-2)n_{0}^2.
\label{CostBCProj_marg0}
\end{equation}
The number of floating-point operations when projecting onto a non-barycentric marginal (as in \eqref{ProjSimplified_BC_marggen}) is
\begin{equation}
    Js(n_{0}n + n_{0}) + \left(3n_{0}+n+2n_{0}n\right) + (2s-2)n_{0}^2.
\label{CostBCProj_marggen}
\end{equation}

\subsection{Projections for Series-Parallel Graph-structured MSBP}\label{subsec:SPProjections}
The series-parallel graph-structured MSBP \eqref{DiscreteSPMSBP} formulated in Sec. \ref{subsec:SeriesParallelGraphStructuredMSBP} involves a series-parallel graph which is not a tree. Recall that the corresponding index set $\Lambda=\Lambda_{\rm{SP}}$ as in \eqref{SPIndexSet}. For any fixed $j\in\llbracket J\rrbracket$, let
$$\Lambda_{\rm{SP}}^{j} := \{j\}\times\left(\llbracket s\rrbracket\setminus\{1,s\}\right).$$
Then, the cost tensor $\bm{C}$ for series-parallel graph-structured MSBP takes the form 
\begin{align}
    \left[\!\bm{C}_{(i_{(j,\sigma)} \mid (j,\sigma)\in\Lambda_{\rm{SP}})}\!\right] = \!\sum_{j=1}^{J}\!\left[\!\bm{C}_{i_{(1,1)},(i_{(j,\sigma)}\mid(j,\sigma)\in\Lambda_{\rm{SP}}^{j}),i_{(1,s)}}^{j}\!\right]   
\label{eqC_SP}
\end{align}
where each $\bm{C}^j\in\left(\mathbb{R}^{n}\right)^{\otimes s}_{> 0}$ is the path-structured cost tensor along the $j$th CPU's path, i.e., 
\begin{align}
    \left[\bm{C}_{i_{(1,1)},(i_{(j,\sigma)}\mid (j,\sigma)\in\Lambda_{\rm{SP}}^{j}),i_{(1,s)}}^{j}\right] = C^{j,1}_{i_{(1,1)},i_{(j,2)}} \nonumber &\\
    + \sum_{\sigma=2}^{s-2}C^{j,\sigma}_{i_{(j,\sigma)},i_{(j,\sigma+1)}} + C^{j,s-1}_{i_{(j,s-1)},i_{(1,s)}}, &
\label{eqCj_SP}
\end{align}
wherein the matrices $C^{j,1}\in\R^{{n}\times{n}}_{\geq 0}$ are the ground cost matrices between marginals $\bm{\xi}^1(\tau_1)$ and $\bm{\xi}^j(\tau_{2})$ for $j\in\llbracket J\rrbracket$, and similarly $C^{j,s-1}$ maps between $\bm{\xi}^j(\tau_{s-1})$ and $\bm{\xi}^1(\tau_{s})$. When $\sigma\in\llbracket s-2\rrbracket\setminus\{1\}$, the matrices $C^{j,\sigma}$ map between $\bm{\xi}^j(\tau_{\sigma})$ and $\bm{\xi}^j(\tau_{\sigma+1})$. 

Consequently, we write $\bm{K}$ in \eqref{eqstructKtensor} as
\begin{subequations}
\begin{align}
    &\!\!\left[\!\bm{K}_{(i_{(j,\sigma)} \mid (j,\sigma)\in\Lambda_{\rm{SP}})}\!\right] \!\!= \!\!\prod_{j=1}^{J}\!\left[\!\bm{K}_{i_{(1,1)},(i_{(j,\sigma)}\mid (j,\sigma)\in\Lambda_{\rm{SP}}^{j}),i_{(1,s)}}^{j}\!\right], \\
    &\left[\bm{K}_{i_{(1,1)},(i_{(j,\sigma)}\mid (j,\sigma)\in\Lambda_{\rm{SP}}^{j}),i_{(1,s)}}^{j}\right] = K^{j,1}_{i_{(1,1)},i_{(j,2)}}, \nonumber \\    &\quad\quad\quad\quad\quad\quad\cdot\left(\prod_{\sigma=2}^{s-2}K^{j,\sigma}_{i_{(j,\sigma)},i_{(j,\sigma+1)}}\right)K^{j,s-1}_{i_{(j,s-1)},i_{(1,s)}},
\end{align}
\label{eqK_SP}
\end{subequations}
and the associated projections can be expressed as follows (proof in Appendix \ref{AppProofPropSPSinkPathCost}).

\begin{proposition}\label{PropSPSinkPathCost}
If $\bm{C}$ has the form \eqref{eqC_SP}, $K^{j,\sigma}:=\exp(-C^{j,\sigma}/\varepsilon)$, $\bm{K}$ as in \eqref{eqK_SP}, and $\bm{U}$ as in \eqref{eqstructUtensor}, then letting
\begin{align}
    A_k := K^{k,1}\left(\prod_{\sigma=2}^{s-1}\diag(\bm{u}_{\sigma}^{k})K^{k,\sigma}\right),\: B_j:=\bigodot_{k\neq j}A_k\:,
\label{eqAk}
\end{align}

and
\begin{align}
    X_\sigma^j &:= K^{j,1}\bigg(\prod_{m=2}^{\sigma-1}\diag(\bm{u}_{m}^j)K^{j,m}\bigg)\:, \nonumber\\
    Y_\sigma^j &:= K^{j,\sigma}\bigg(\prod_{m=\sigma+1}^{s-1}\diag(\bm{u}_{m}^j)K^{j,m}\bigg)\:, \nonumber\\
    Z_{\sigma_1,\sigma_2}^j &:= K^{j,\sigma_1}\bigg(\prod_{m=\sigma_1+1}^{\sigma_2-1}\diag(\bm{u}_m^j)K^{j,m}\bigg)\:, \nonumber
\end{align}
(so $A_k=X_s^k=Y_1^k$),
the projection \eqref{DefprojsigmaComponent} takes the form
\begin{subequations}
\begin{align}   {\rm{proj}}_{(1,1)}\!\left(\bm{K}\odot\bm{U}\right)\! &= \bm{u}_1^1\odot\left(\bigodot_{k=1}^JA_k\right)\bm{u}_s^1, \label{ProjSimplified_SP_marg1} \\    {\rm{proj}}_{(1,s)}\!\left(\bm{K}\odot\bm{U}\right)\! &= \left({\bm{u}_1^1}^{\top}\cdot\bigodot_{k=1}^JA_k\right)^{\top}\odot\bm{u}_s^1, \label{ProjSimplified_SP_margs} \\   {\rm{proj}}_{(j,\sigma)}\!\left(\bm{K}\odot\bm{U}\right)\! &= \bm{u}_\sigma^j\!\odot\!  \diag\Big(\!Y_\sigma^j\diag(\bm{u}_s^1)B_j^\top\diag({\bm{u}_1^1})X_\sigma^j\!\Big), \label{ProjSimplified_SP_marggen}
\end{align}
\label{ProjSimplified_SP}
\end{subequations}
where $(j,\sigma)\in\bigcup_{k\in\llbracket J\rrbracket}\Lambda_{\rm{SP}}^k$ in \eqref{ProjSimplified_SP_marggen}. Furthermore, the projection \eqref{DefprojsigmaComponentDouble} can be expressed as
\begingroup\allowdisplaybreaks
\begin{subequations}
\begin{align}
    &{\rm{proj}}_{(1,1),(j,2)}\!\left(\bm{K}\odot\bm{U}\right)\! = \diag\left(\bm{u}_1^1\right) \nonumber\\
    &\quad\quad\cdot \Bigg(K^{j,1}\odot\bigg(B_j\diag(\bm{u}_s^1){Y_2^j}^{\top}\bigg)\Bigg) \cdot \diag\left(\bm{u}_{2}^{j} \right), \label{Proj2Simplified_SP_marg1}\\
    &{\rm{proj}}_{(j,s-1),(1,s)}\!\left(\bm{K}\odot\bm{U}\right)\! = \diag\left(\bm{u}_{s-1}^{j}\right) \nonumber\\
    &\quad\cdot \Bigg(K^{j,s-1}\odot\bigg({X_{s-1}^j}^\top\diag(\bm{u}_1^1)B_j\bigg)\Bigg) \cdot \diag\left(\bm{u}_{s}^1\right),\label{Proj2Simplified_SP_margs}\\
    &{\rm{proj}}_{(j,\sigma),(j,\sigma+1)}\!\left(\bm{K}\odot\bm{U}\right)\! = \diag\left(\!\bm{u}_{\sigma}^{j}\!\right) \nonumber\\
    &\!\cdot\! \Bigg(K^{j,\sigma}\odot\bigg({X_{\sigma}^j}^\top\diag(\bm{u}_1^1)B_j\diag(\bm{u}_s^1){Y_{\sigma+1}^j}^{\!\!\!\!\!\!\top}\bigg)\Bigg) \!\cdot\! \diag\left(\!\bm{u}_{\sigma+1}^j\!\right). \label{Proj2Simplified_SP_marggen}
\end{align} 
\label{Proj2Simplified_SP}
\end{subequations}
\endgroup
\end{proposition}

Directly computing $A_k$ in \eqref{eqAk} requires $1+2(s-2)$ matrix-matrix multiplications, each in this case involving $n^3$ floating-point operations. 
Thus, \eqref{ProjSimplified_SP_marg1} and \eqref{ProjSimplified_SP_margs} each involve
\begin{equation}
    J(1+2(s-2))n^3+(J+1)n^2+n
\label{CostSPProj_marg0s}
\end{equation}
operations, and \eqref{ProjSimplified_SP_marggen} involves
\begin{equation}
    \left(J(1+2(s-2))+3\right)n^3+(J-1)n^2+n
\label{CostSPProj_marggen}
\end{equation}
operations. 

\begin{remark}\label{ComplexityComparison}
For $n_{0} = n, s \geq 2, J \geq 2$, we note that 
$$\text{the expression}\,\eqref{CostBCProj_marggen} = (Js + 2s) n^2 + (Js + 4) n \geq 8n^2 + 8n.$$
In contrast, the expression $\eqref{CostSPProj_marggen} = (2Js - 3J + 3)n^3 + (J-1)n^2+n \!\geq\! n^3 \!+\! n^2 \!+\! n$. 
Thereby for $n$ large,
Propositions \ref{PropBCSinkPathCost} and \ref{PropSPSinkPathCost} enable computation of unimarginal projections with $\mathcal{O}(Jsn^2)$ complexity for the barycentric case, and $\mathcal{O}(Jsn^3)$ for the series-parallel case. The $n^3$ complexity of the latter is a consequence of the treewidth ($=2$) of the series-parallel graph structure. Both cases have complexity linear in $Js$.
\end{remark}


\subsection{Overall Algorithm}\label{ss:Alg.OA}
Having shown that the computation of the projections are tractable in the number of marginals, we are able to efficiently solve \eqref{DiscreteMSBP} for all of our proposed graph structures: path, barycentric (BC), and series-parallel (SP). In applying MSBP for computational resource prediction, the choice of problem graph structure can affect the quality of the prediction. For single-core software, all graph structures necessarily degenerate to the path structure. But for a multi-core software, there are two structurally distinct options: BC and SP. In the following Sec. \ref{sec:Experiements}, we evaluate the merits of these structures.

Regardless of the choice of graph structure, our overall methodology is as follows.

\noindent\textbf{Step 1.} Execute the software of interest $n$ times over $[0,t]$, generating hardware resource state snapshots $\{\boldsymbol{\xi}^{i,j}\left(\tau_\sigma\right)\}_{i=1}^{n}$ for all $(j,\sigma)\in\llbracket J\rrbracket\times\llbracket s\rrbracket$. Our marginal distributions $\mu_\sigma^j$ are then as per \eqref{DefEmpiricalMeasure}. \\
\noindent\textbf{Step 2.} Pick an appropriate graph structure (path, BC, or SP). From the marginals in \textbf{Step 1}, construct the bimarginal Euclidean distance matrices $C^{j,\sigma}$. If the BC structure is used, compute the barycenters $\mu_\sigma^0$. Run the Sinkhorn recursion \eqref{MultimarginalSink} until convergence, i.e. until the Hilbert projective metric \eqref{HilbertMetricPosOrthant} between the new and old $\bm{u}_\sigma^j$ is less than some numerical tolerance for all marginals. The resulting $\bm{M}_{\rm{opt}}$ given by \eqref{Mopt} is the solution for our MSBP. \\
\noindent\textbf{Step 3.} Given any time $\tau\in[0,t)$ and CPU index $j$, use \eqref{MuInterpolation} to predict our estimate $\widehat{\mu}_{\tau}^j$ for the hardware resource usage distribution of CPU $j$ at time $\tau$.


\section{Experiments}\label{sec:Experiements}
In Sec. \ref{sec:GraphStructuredMSBP} and \ref{sec:Algorithms}, we explained how to formulate and efficiently solve graph-structured MSBPs to compute the optimal mass transport plan $\bm{M}^{\rm{opt}}$ between marginal distributions \eqref{DefEmpiricalMeasure} which in our application correspond to the statistics of computational resource usage at times \eqref{SnapshotTimes}. 

In Remark \ref{RemarkInterpolation}, we outlined how the computed $\bm{M}^{\rm{opt}}$ can be used to make predictions on a software's resource usage at any (possibly out-of-sample) given time. To evaluate the quality of these predictions, in this Section we perform two experiments: one with a single-core software and one with a multi-core software. We quantify the quality of our predictions against known resource usage data.

For all experiments herein, benchmark programs were profiled on an Intel Xeon E5-2683 v4 processor with 16 cores and a 40MB LLC running Ubuntu 16.04.7. We leveraged Intel's Cache Allocation Technology (CAT) \cite{intel-2015-cat} and Memguard \cite{Yun16-memguard-journal} to control the amount of LLC and memory bandwidth available to the benchmark programs.

The implementation of the algorithms described in Sec. \ref{sec:Algorithms}, the solution of the MSBP \eqref{GenericDiscreteMSBP} itself, and all other processing of profiled data, were done in MATLAB R2019b on a Debian 12 Linux machine with an AMD Ryzen 7 5800X CPU. 

\vspace{-1mm}
\subsection{Single Core Experiment}\label{subsec:SingleCoreExperiment}
For this experiment we used a custom single-core software\footnote{Github repository: {\tiny{\url{https://github.com/abhishekhalder/CPS-Frontier-Task3-Collaboration}}}} which we wrote in C language, implementing a path-following nonlinear model predictive controller (NMPC) for a kinematic bicycle model (KBM) \cite{kong2015kinematic,haddad2022density}. At each control step, the NMPC used the IPOPT nonlinear program solver \cite{IPOPT} to determine the control decision which minimizes the sum of various measurements of deviation from the desired path while promoting smoothness of the control inputs. For details on this control software, we refer the readers to \cite[Sec. V]{bondar2023path}. 

\subsubsection{Profiling} To gather the execution profiles for our NMPC control software, we ran our application on an isolated CPU and used the Linux perf tool \cite{perf}, version 4.9.3, to sample the computational resource state vector $\bm{\xi}$ as in \eqref{defFeatureVec} every 10 ms. 

For each profile, we ran the NMPC software for $n_c:=5$ control steps with a fixed reference trajectory and a fixed combination of LLC and memory bandwidth partitions (30 MB and 30 MBps respectively, both allocated in blocks of 2 MB) to ensure a fixed execution context. 

\begin{figure*}[t]
    \centering
\includegraphics[width=0.95\linewidth]{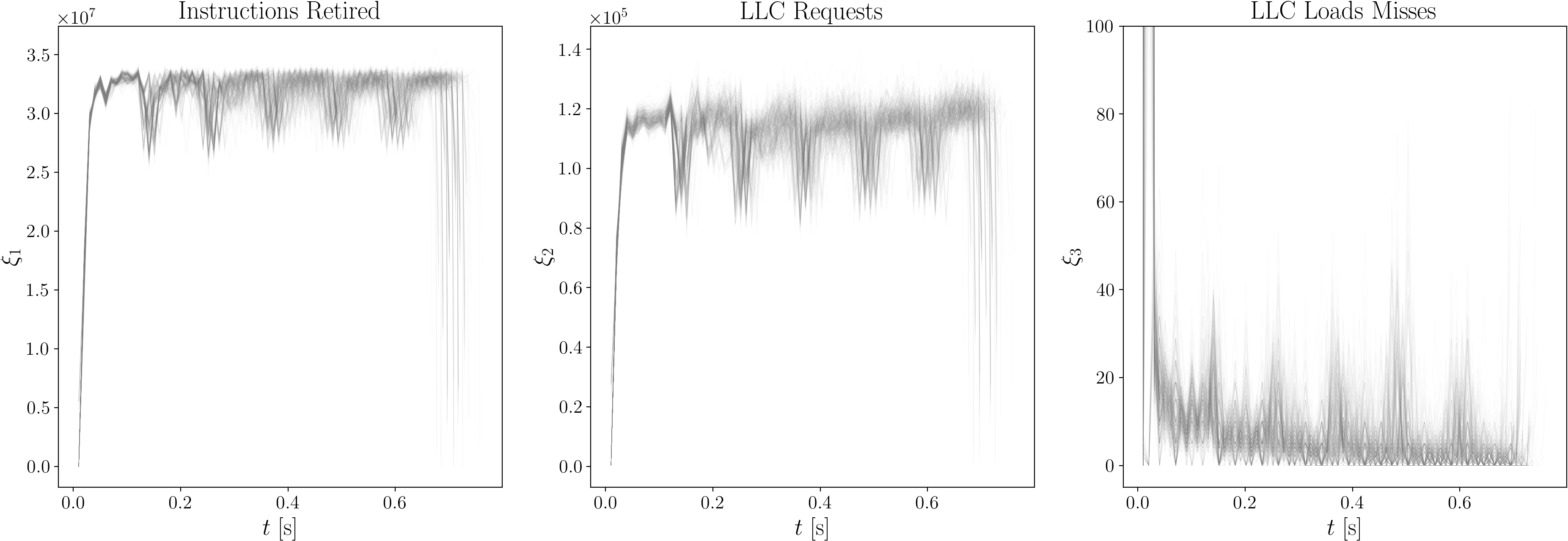}    \caption{{\small{Components of the measured feature vector $\bm{\xi}$ in \eqref{defFeatureVec} for the single core experiment in Sec. \ref{subsec:SingleCoreExperiment}, for all of the five control cycles for 500 executions of the NMPC software, for a fixed cyber-physical context. 
}}}
\vspace*{-0.15in}
\label{FigAllMeasuredStatesForOneContext}
\end{figure*}

For our chosen context, we profiled the controller $n:=500$ times. Fig. \ref{FigAllMeasuredStatesForOneContext} overlays all of the profiles, split into its three components as per \eqref{defFeatureVec}. \\

\subsubsection{Numerical Results}
To account for the asynchrony among the profiles, we analyzed the statistics of the control cycle end times, i.e., the (random) times needed to complete one pass of the NMPC feedback loop. For details, we refer the readers to \cite[Table I, Fig. 4]{bondar2023path}. We then consider snapshots of the profiling data at times corresponding to the average control cycle end times as well as at $s_{\rm{int}}$ equi-spaced intervals between control cycles. This results in $s=1+n_c(s_{\rm{int}}+1)$ total marginals forming our path-structured MSBP \eqref{DiscreteMSBP}, and completes {\bf{Step 1}} in section \ref{ss:Alg.OA}.

\begin{figure}[t]
    \centering    \includegraphics[width=0.9\linewidth]{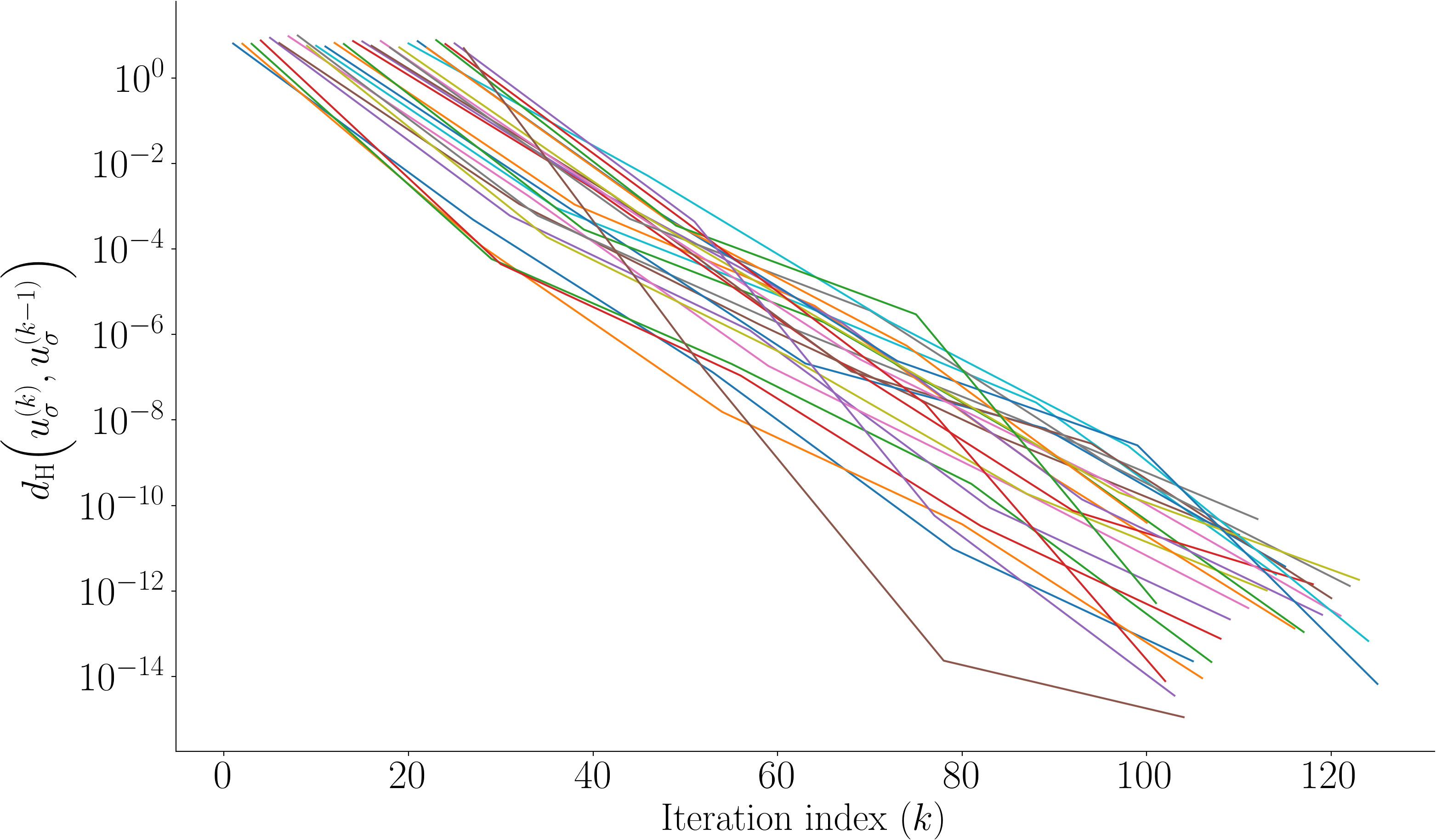}
    \caption{{\small{Linear convergence of Sinkhorn iterations \eqref{MultimarginalSink} for the single core experiment in Sec. \ref{subsec:SingleCoreExperiment}, for $s_{\rm{int}}=4$ w.r.t. the Hilbert's projective metric $d_{\rm{H}}$ in \eqref{HilbertMetricPosOrthant} between $\bm{u}_{\sigma\in\llbracket s\rrbracket}$ at iteration indices $k$ and $k-1$.}}}
\vspace*{-0.15in}
\label{FigConvergenceHilbert}
\end{figure}

Setting the regularizer $\varepsilon=0.1$ and a numerical tolerance of $10^{-14}$, we solve the discrete path-structured MSBP \eqref{DiscreteMSBP} with squared Euclidean cost $\bm{C}$ as in \eqref{eqC_SingleCore}. Fig. \ref{FigConvergenceHilbert} depicts the linear convergence of the Sinkhorn iterations (as in Proposition \ref{propSinkhornConvergence}) shown w.r.t. the Hilbert metric \eqref{HilbertMetricPosOrthant}. The path structure of the problem yields efficient computation of the Sinkhorn iterations leading to rapid convergence. In particular, setting $s_{\rm{int}}:=4$ (and so $s=26$), we solve the MSBP \eqref{DiscreteMSBP} with $n^{s} = 500^{26}$ decision variables in approx. 10s and $\approx 120$ Sinkhorn iterations with no optimizations to the code. Notice that problems of this size are impractical to load in memory, let alone to solve efficiently using off-the-shelf solvers. This completes {\bf{Step 2}} in section \ref{ss:Alg.OA}.

\begin{figure}[h]
    \centering    \includegraphics[width=0.9\linewidth]{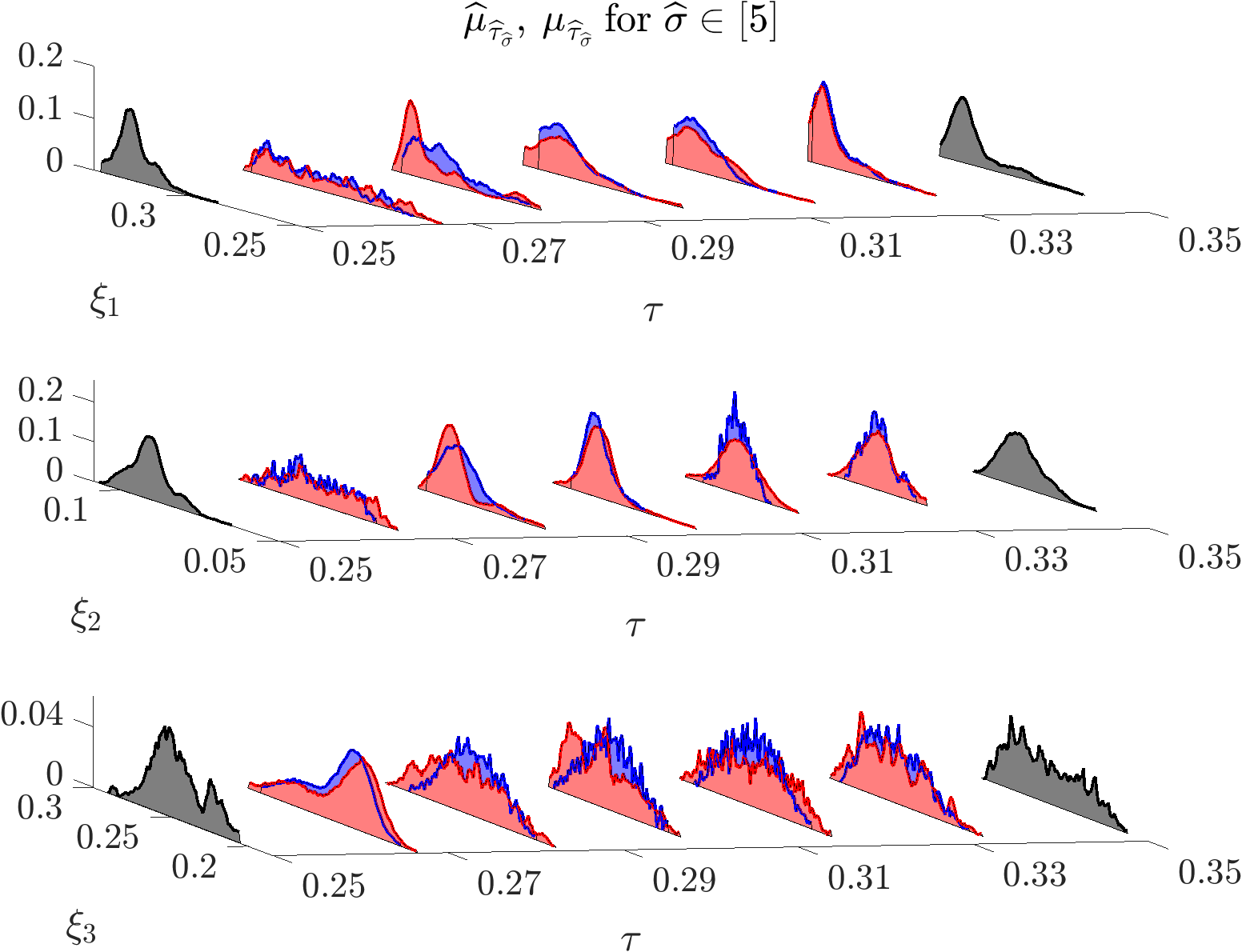}
    \caption{{\small{Predicted $\widehat{\mu}_{\widehat{\tau}_{\widehat{\sigma}}}$ (\emph{blue}) vs. measured $\mu_{\widehat{\tau}_{\widehat{\sigma}}}$ (\emph{red}) at times $\hat{\tau}_{\widehat{\sigma}\in\llbracket 5\rrbracket}$ for the single core experiment in Sec. \ref{subsec:SingleCoreExperiment}, during the 3rd control cycle with $s_{\rm{int}}=4$. Distributions at the control cycle boundaries are in \emph{black}.}}}
\vspace*{-0.15in}
\label{FigInterpMarginals1D}
\end{figure}

For {\bf{Step 3}}, we compute $s_{\rm{int}}+1=5$ distributional predictions $\widehat{\mu}_{\widehat{\tau}_{\widehat{\sigma}}}$ at times $\widehat{\tau}_{\widehat{\sigma}}$, temporally equispaced throughout the duration of the 3rd control cycle, i.e., between $\tau_{2(s_{\rm{int}}+1)+1}$ and $\tau_{3(s_{\rm{int}}+1)+1}$, with
$$ \widehat{\tau}_{\widehat{\sigma}}=\tau_{2(s_{\rm{int}}+1)+1} + \Bigg(\frac{\tau_{3(s_{\rm{int}}+1)+1}-\tau_{2(s_{\rm{int}}+1)+1}}{s_{\rm{int}}+2}\Bigg)\widehat{\sigma}, $$
where the index $\widehat{\sigma}\in\llbracket s_{\rm{int}}+1 \rrbracket$. Since $\widehat{\tau}_{\widehat{\sigma}}\in[\tau_{2(s_{\rm{int}}+1)+\widehat{\sigma}},\tau_{2(s_{\rm{int}}+1)+\widehat{\sigma}+1}]$, we used \eqref{MuInterpolation} with $\sigma=2(s_{\rm{int}}+1)+\widehat{\sigma}$ to arrive at the predictions. In Fig. \ref{FigInterpMarginals1D}, we compare our predictions against the observed empirical distributions.

Fig. \ref{FigInterpMarginals1D} shows that our predictions capture the mode(s) of the measured distributions quite well. Further improvements can be made by placing marginals closer together in time. In this experiment, we do this by increasing the value of $s_{\rm{int}}$, i.e., by placing more marginals equi-spaced between control cycle boundaries. In Table \ref{tab:WassersteinErrorsVaryingSint}, we report the Wasserstein distances $W(\cdot,\cdot)$ as in \eqref{DefWasserstein} between the corresponding predicted and measured distributions:
\begin{align}
W_{\widehat{\sigma}} := W(\widehat{\mu}_{\widehat{\tau}_{\widehat{\sigma}}},\mu_{\widehat{\tau}_{\widehat{\sigma}}}) \quad\forall \widehat{\sigma}\in\llbracket s_{\rm{int}}+1 \rrbracket.
\label{Wasserstein}
\end{align}
We computed \eqref{Wasserstein} as the square roots of the optimal values of the corresponding Kantorovich linear programs \cite[Ch. 3.1]{peyre2019computational} that results from specializing \eqref{DiscreteMSBP} with $s=2$, $\varepsilon=0$. As expected, placing the marginals closer together in time results in higher accuracy in predictions.

\begin{table}[t]
\centering
\scriptsize
\begin{tabulary}{\textwidth}{|c|c|c|c|c|c|} 
 \hline
 $s_{\rm{int}}$ & $W_{1}$ & $W_{2}$ & $W_{3}$ & $W_{4}$ & $W_{5}$\\
 \hline\hline
 $0$ & $2.049\times 10^{-4}$ & - & - & - & - \\
 \hline
 $1$ & $2.270\times10^{-4}$ & $1.175\times10^{-4}$ & - & - & - \\
 \hline
 $2$ & $5.772\times10^{-4}$ & $9.163\times10^{-5}$ & $3.794\times10^{-5}$ & - & - \\
 \hline
 $3$ & $2.241\times10^{-4}$ & $1.643\times10^{-4}$ & $1.234\times10^{-4}$ & $6.010\times10^{-5}$ & -  \\
\hline
 $4$ & $6.372\times10^{-5}$ & $1.2691\times10^{-4}$ & $9.176\times10^{-5}$ & $6.689\times10^{-5}$ & $2.111\times10^{-5}$ \\
\hline
\end{tabulary}
\caption{{\small{For the single core experiment in Sec. \ref{subsec:SingleCoreExperiment}, the number of intracycle marginals $s_{\rm{int}}$ vs. Wasserstein distances $W_{\widehat{\sigma}}$ as in \eqref{Wasserstein}. 
}}}
\label{tab:WassersteinErrorsVaryingSint}
\vspace*{-0.15in}
\end{table}

\subsection{Multi-core Experiment}\label{subsec:MultiCoreExperiment}
\begin{figure}[t]
    \centering    \includegraphics[width=0.98\linewidth]{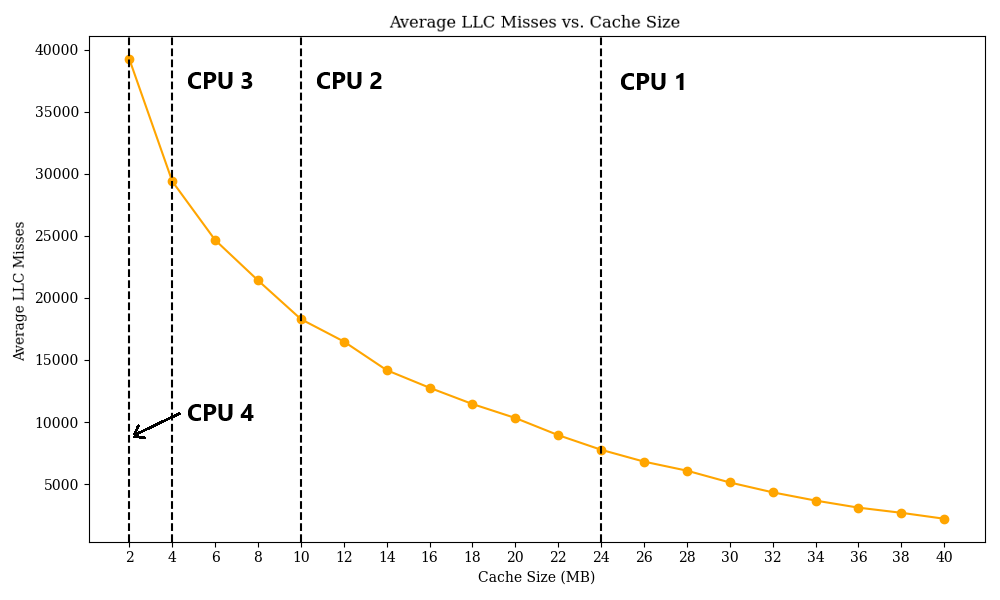}
    \caption{{\small{LLC miss rate per cache size for the Canneal benchmark in Sec. \ref{subsec:SingleCoreExperiment}. The dashed vertical lines show our resource allocations per core.}}}
\vspace*{-0.15in}
\label{FigCannealMiss}
\end{figure}

For this experiment, we used the Canneal benchmark from the PARSEC suite~\cite{bienia2008parsec}. Canneal is a resource-heavy, multi-threaded application that simulates an anneal workload to minimize routing costs for chip design. Specifically, Canneal pseudorandomly picks pairs of input elements to swap in a tight loop, leading to a heavy reliance on both cache and memory bandwidth. This can be seen in Fig.~\ref{FigCannealMiss} where the average number of LLC misses decreases as the cache size increases. The dashed vertical lines represent the cache allocation used in our profiling.

\subsubsection{Profiling} To enable multi-core profiling, we made a small modification to Canneal's source code to pin each created thread to its own core. We then use perf, as with the single core experiment, to sample the computational resource state $\bm{\xi}^{j}$, $j\in\llbracket J\rrbracket$, as in \eqref{defFeatureVec} for every 10 ms. 

\begin{figure*}[t]
    \centering    \includegraphics[width=0.99\linewidth]{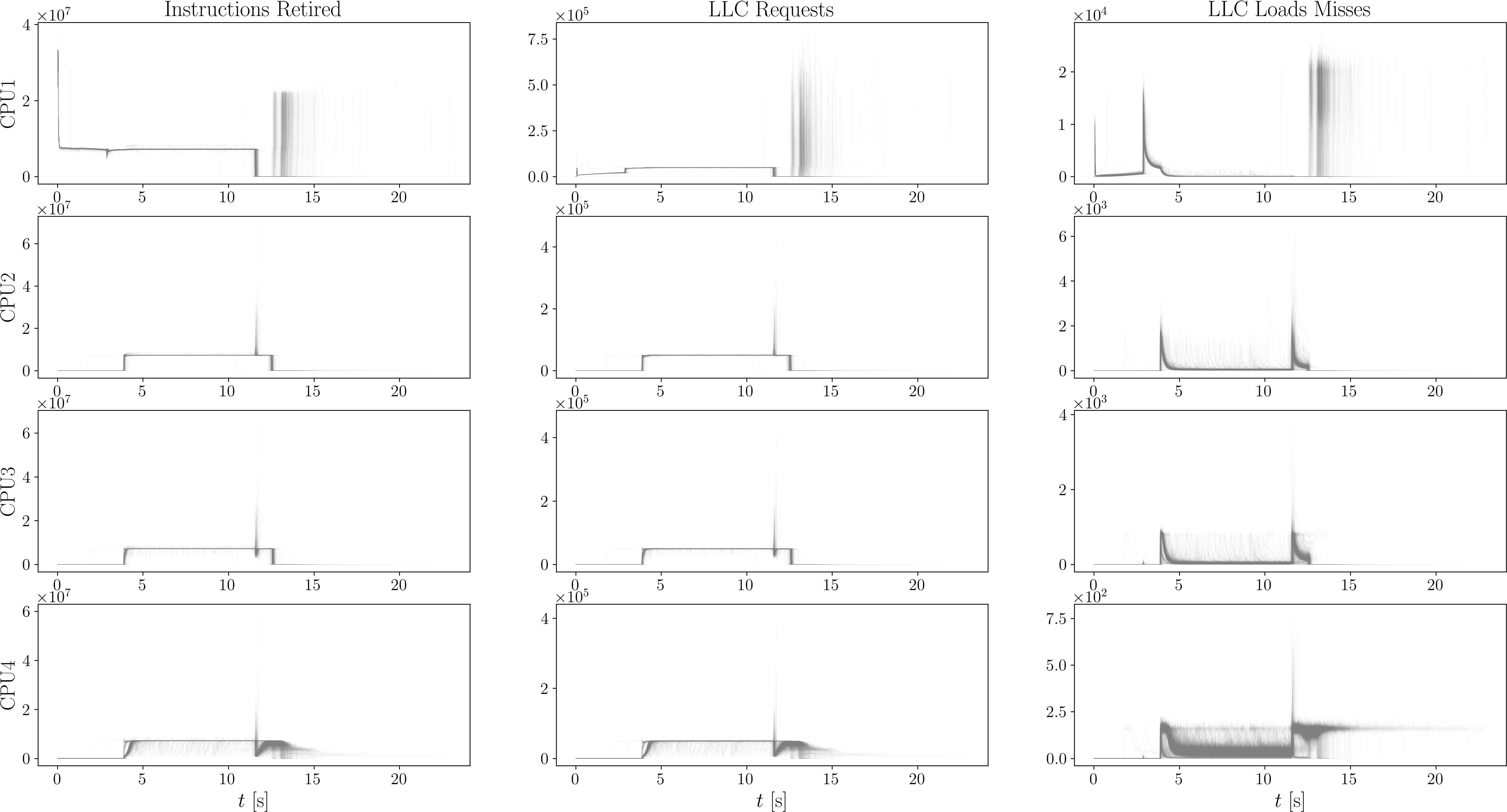}
    \caption{{\small{Components of $\bm{\xi}^{j}$, $j\in\llbracket 4\rrbracket$, for the $n=400$ executions of the Canneal benchmark in Sec. \ref{subsec:MultiCoreExperiment}, for all CPUs. Observe the erratic behavior of CPU4, which has the least hardware resources, as well that of CPU1, which has the most.}}}
\vspace*{-0.15in}
\label{FigAllMeasuredStatesCanneal}
\end{figure*}

Desiring to examine how Canneal behaves when its resource allocation varies across its cores, we profile using a single context defined by $J=4$ cores, and resource allocations of $(24,10,4,2)$ MBs of LLC and $(125,25,5,1)$ MBps of memory bandwidth, ordered left-to-right by increasing CPU number.
So, CPU 1 has the highest resource allocation while CPU 4 has the lowest. 
For this choice of context, we collected a total of $n=400$ profiles for the Canneal software; these profiles are shown in Fig. \ref{FigAllMeasuredStatesCanneal}. For the barycentric formulation, we fix $n_0=600$.

\subsubsection{Numerical Results}

We placed snapshots at times $\tau_\sigma\in\{0.0, 0.5, 1.5, 2.5, 5.0, 9.5, 10.5\}$ (so $s=7$) and interpolated at times $\widehat{\tau}_{\widehat{\sigma}}\in\{0.8, 2.2, 7.0, 9.0, 10.0\}$ (so $\widehat{\sigma}\in\llbracket 5\rrbracket$). Both BC and SP methods were then used to solve the MSBPs \eqref{DiscreteBaryMSBP} and \eqref{DiscreteSPMSBP} respectively, therefrom we used \eqref{MuInterpolation} to estimate the resource usage distributions at the desired times.


Once again, the code for solution of \eqref{DiscreteBaryMSBP} and \eqref{DiscreteSPMSBP}, including the computation of the projections in Propositions \ref{PropBCSinkPathCost} and \ref{PropSPSinkPathCost}, were implemented in MATLAB with some minor optimizations owing to the large size of the problem (the BC and SP formulations have $400^{35}$ and $400^{22}$ decision variables respectively; see Sec. \ref{sec:Algorithms}). 

\begin{figure}[t]
    \centering    \includegraphics[width=0.99\linewidth]{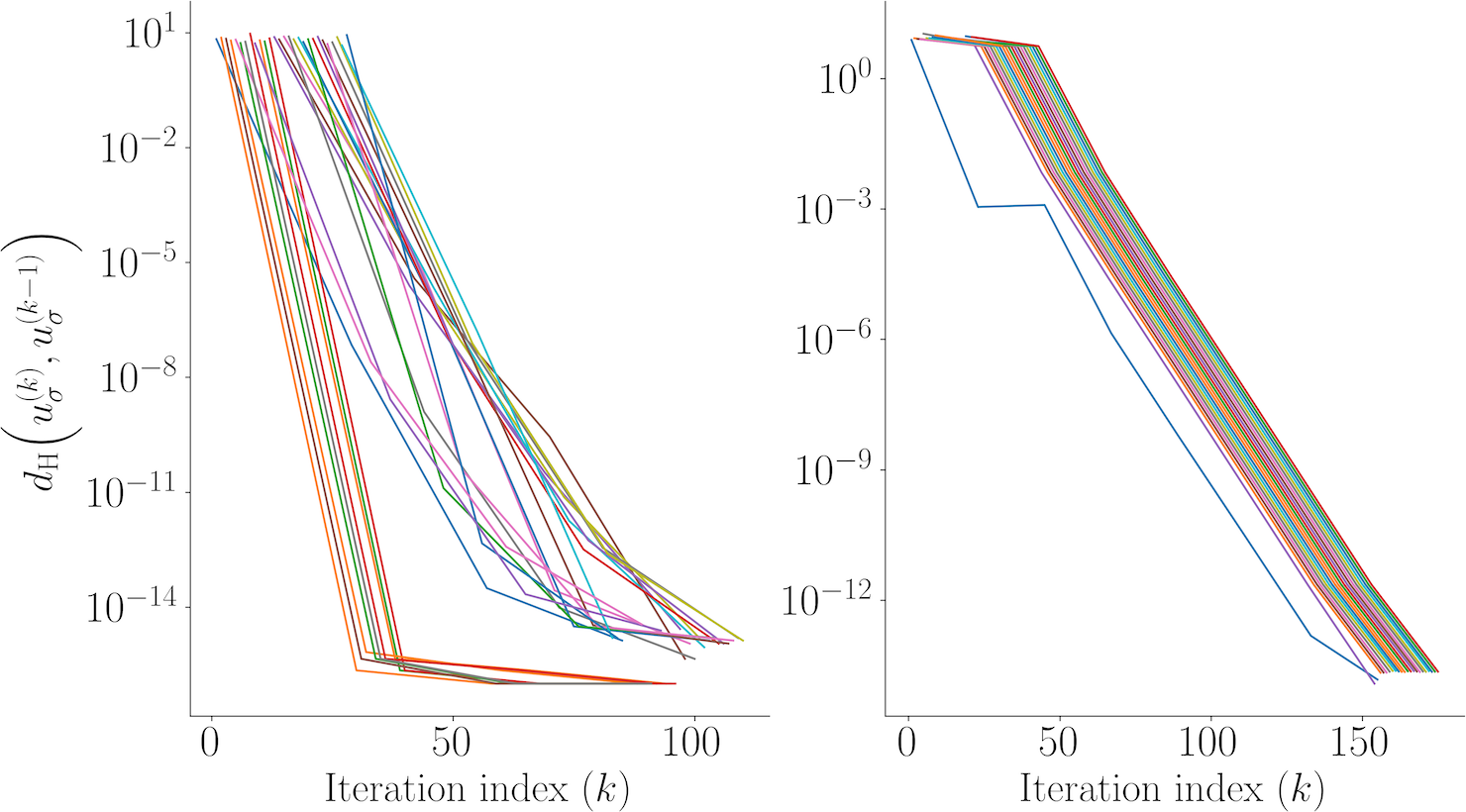}
    \caption{{\small{Convergence of Sinkhorn iterations of both BC (left) and SP (right) algorithms for the Canneal benchmark in Sec. \ref{subsec:MultiCoreExperiment}, shown w.r.t. the Hilbert projective metric $d_{\rm{H}}$ as in \eqref{HilbertMetricPosOrthant}.}}}
\vspace*{-0.15in}
\label{FigConvergenceHilbert_Canneal}
\end{figure}

\begin{figure*}
    \centering
    \begin{subfigure}[b]{0.475\textwidth}
        \centering
        \includegraphics[width=\textwidth]{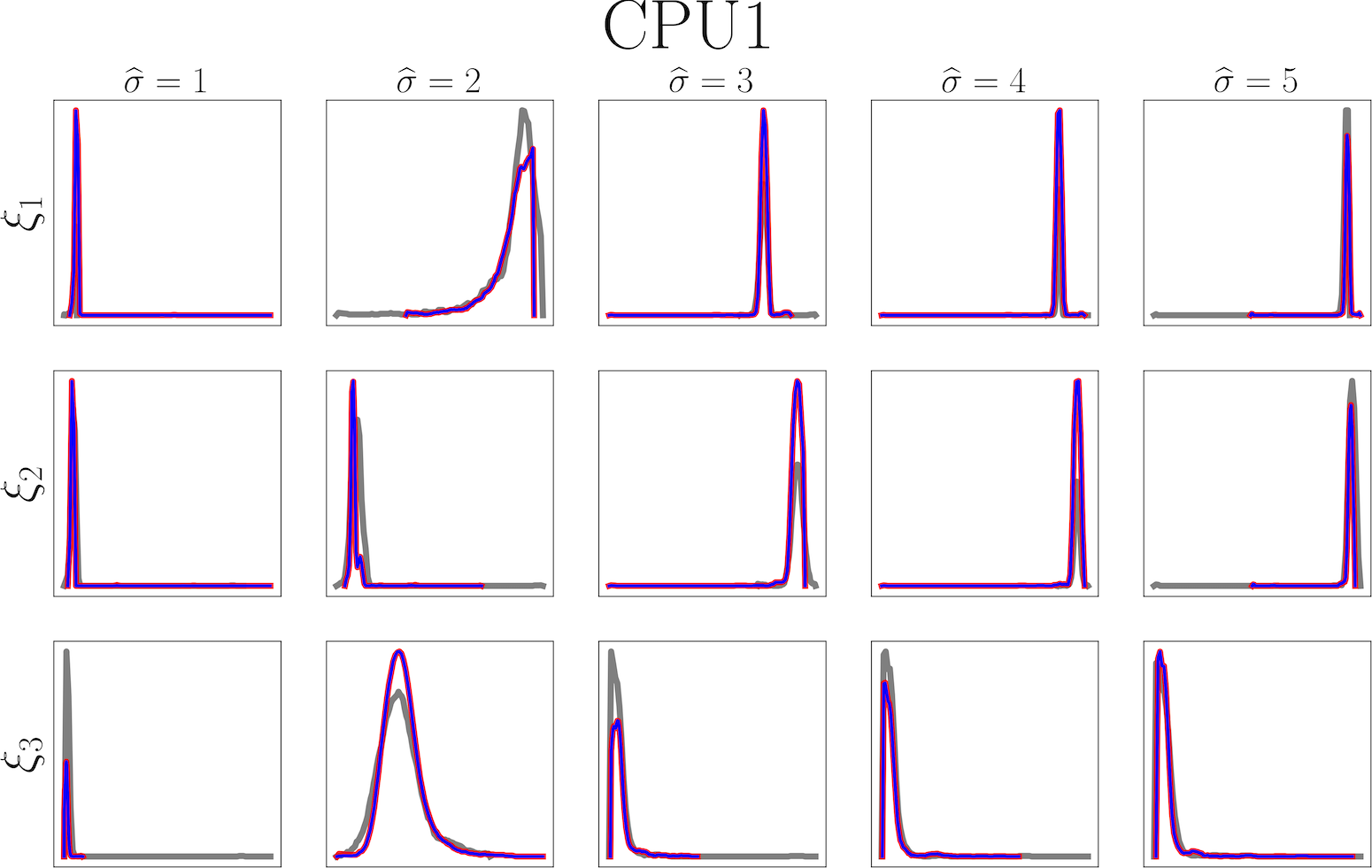}
    \end{subfigure}
    \hfill
    \begin{subfigure}[b]{0.475\textwidth}  
        \centering 
        \includegraphics[width=\textwidth]{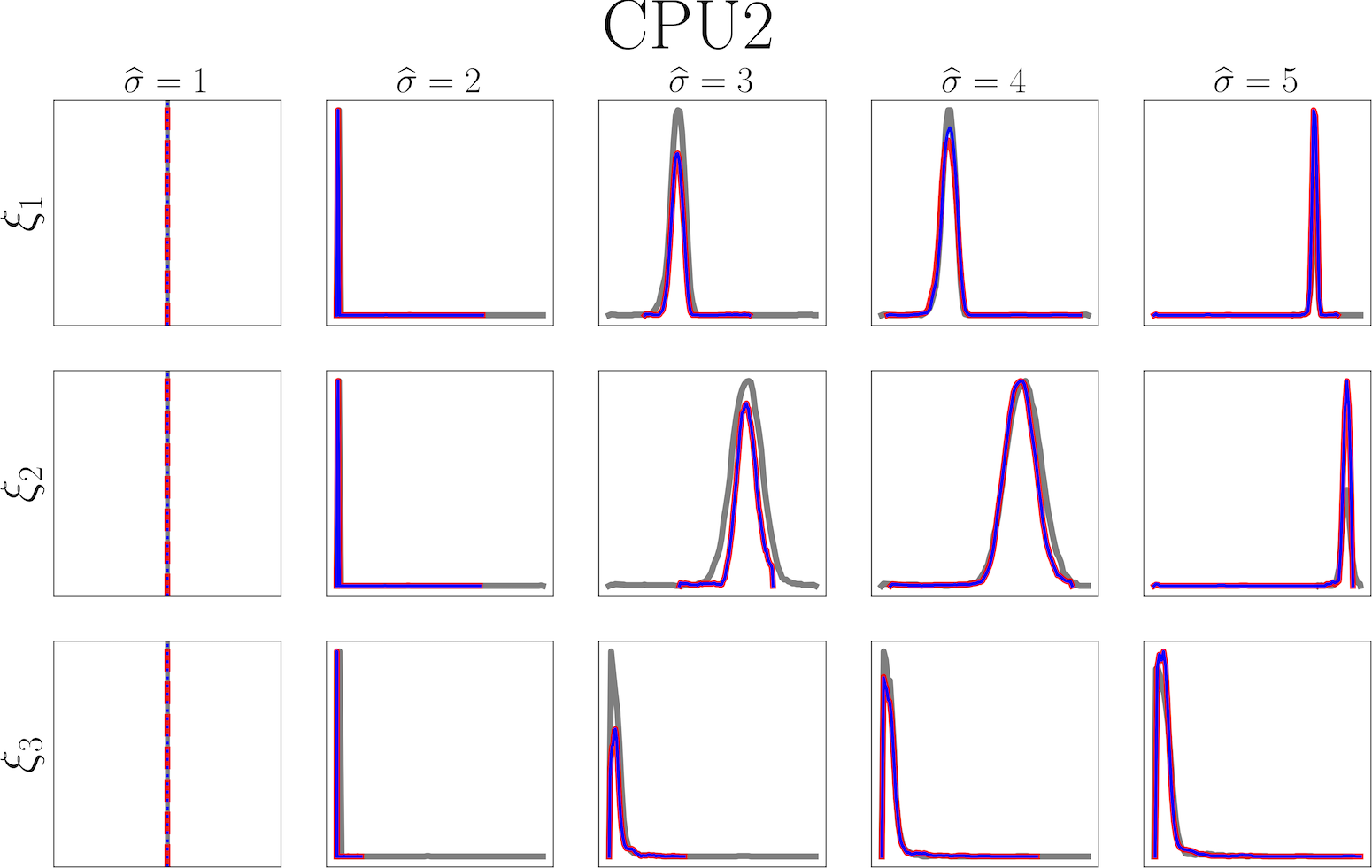}
    \end{subfigure}
    \vskip\baselineskip
    \begin{subfigure}[b]{0.475\textwidth}   
        \centering 
        \includegraphics[width=\textwidth]{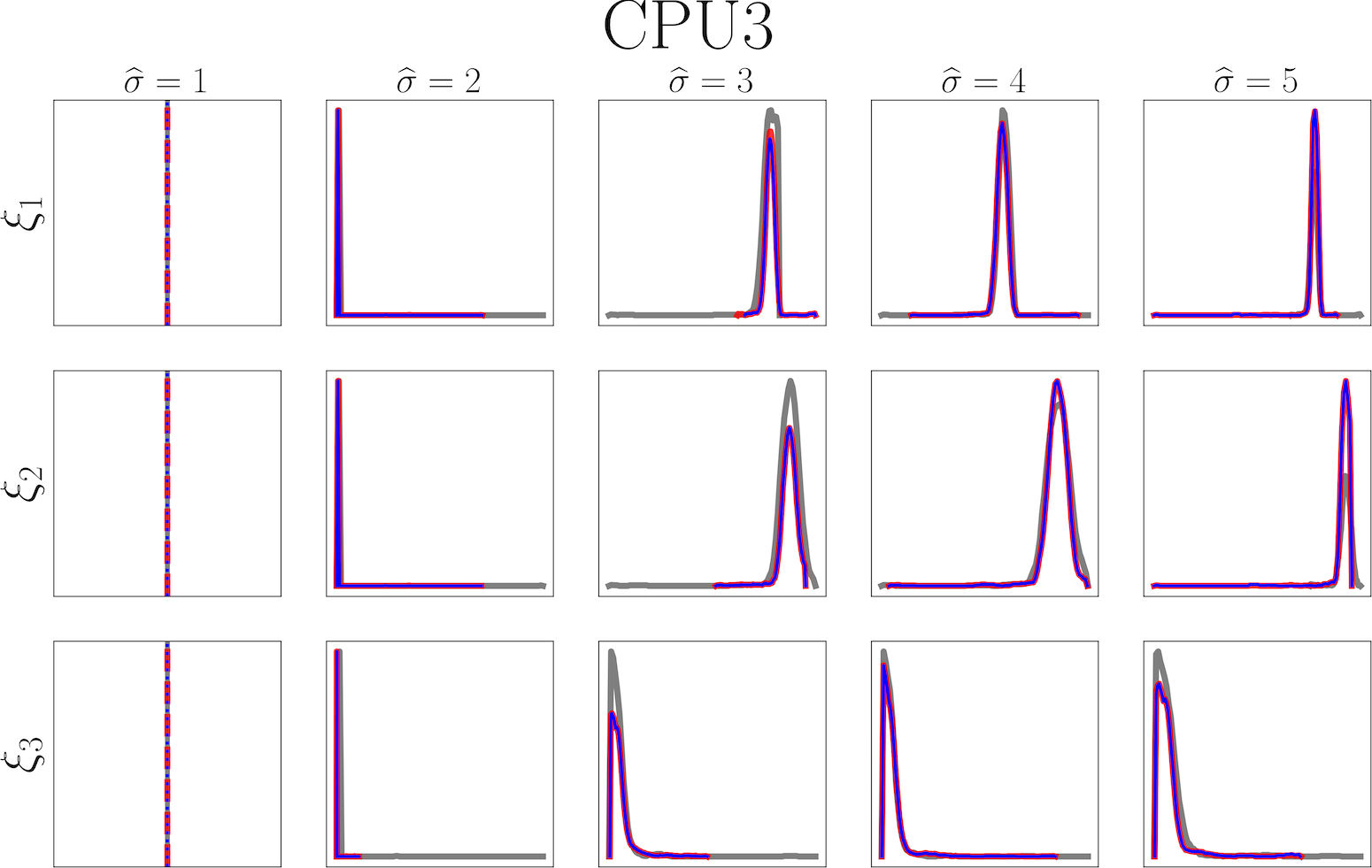}
    \end{subfigure}
    \hfill
    \begin{subfigure}[b]{0.475\textwidth}   
        \centering 
        \includegraphics[width=\textwidth]{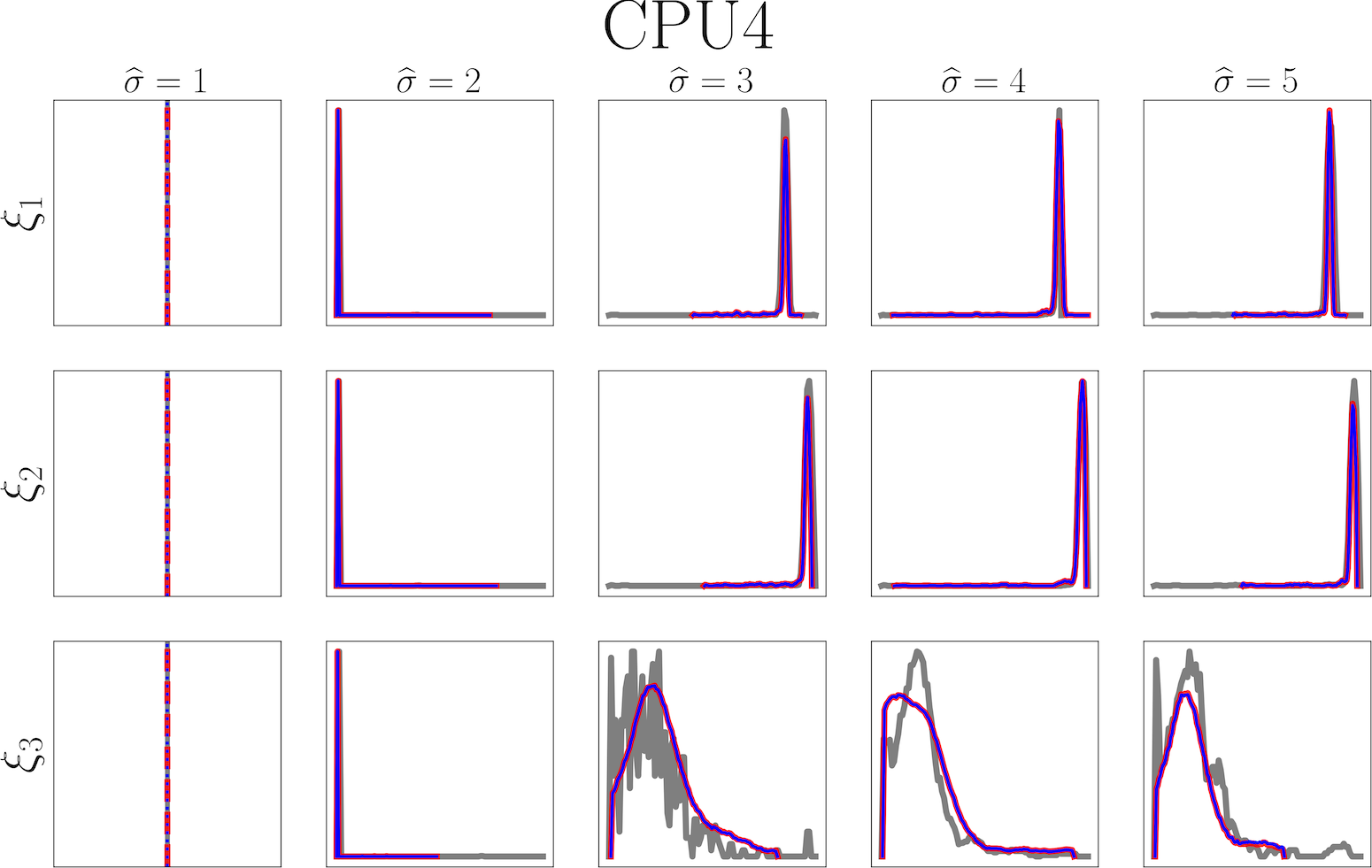}
    \end{subfigure}
    \caption{{\small{Resource usage distributions predicted with the BC structure (red), the SP structure (blue), vs. the measured distributions (grey) at prediction times $\widehat{\tau}_{\widehat{\sigma}}\in\llbracket 5\rrbracket$ for the multi-core Canneal benchmark in Sec. \ref{subsec:MultiCoreExperiment}.}}}
    \label{FigInterpMarginals_Canneal}
    \vspace*{-0.15in}
\end{figure*}

With the regularizer $\varepsilon=0.05$ and a numerical tolerance of $10^{-13}$, the BC and SP algorithms converged in $0.4810$s and $0.5055$s, with $110$ and $127$ iterations respectively (see Fig. \ref{FigConvergenceHilbert_Canneal}). The predicted distributions, marginalized to each component of $\bm{\xi}^{j}$, $j\in\llbracket J\rrbracket$, are shown in Fig. \ref{FigInterpMarginals_Canneal}. Finally, Table \ref{tab:WassersteinErrorsBC_SP} compares the Wasserstein distances of the predicted distributions for each CPU from those measured from our profiles, i.e.,
\begin{align}
W_{\widehat{\sigma}}^j := W(\widehat{\mu}_{\widehat{\tau}_{\widehat{\sigma}}}^j,\mu_{\widehat{\tau}_{\widehat{\sigma}}}^j) \quad\forall (j,\widehat{\sigma})\in\llbracket J\rrbracket\times\llbracket s_{\rm{int}}+1 \rrbracket.
\label{WassersteinBC_SP}
\end{align}

These figures show the same behavior as in the single-core experiment -- primarily, that our predictions accurately capture the mode(s) of the resource usage distribution at all interpolation times. Note that at $\widehat{\sigma}=1,\:\widehat{\tau}_{\widehat{\sigma}}=0.8$, all CPUs except the CPU 1 are idle. Our predictions match the measured distribution exactly for the CPUs 2, 3, and 4 (i.e. the Dirac delta at $0$ for all components), but not so for CPU 1, which has a nontrivial distribution. 

\begin{table}[t]
\centering
\scriptsize
\begin{tabulary}{\textwidth}{|c|c|c|c|c|c|} 
 \hline
 $j$ & $W_{1}^j$ & $W_{2}^j$ & $W_{3}^j$ & $W_{4}^j$ & $W_{5}^j$\\
 \hline\hline
 $1$ & $4.077\times10^{-5}$ & $1.009\times10^{-7}$ & $2.131\times10^{-7}$ & $1.976\times10^{-7}$ & $1.509\times10^{-7}$ \\
 \hline
 $2$ & $0$ & $1.135\times10^{-7}$ & $2.342\times10^{-7}$ & $7.684\times10^{-8}$ & $8.805\times10^{-8}$ \\
 \hline
 $3$ & $0$ & $1.149\times10^{-7}$ & $1.534\times10^{-7}$ & $5.752\times10^{-8}$ & $6.538\times10^{-8}$  \\
\hline
 $4$ & $0$ & $3.647\times10^{-8}$ & $2.146\times10^{-7}$ & $1.906\times10^{-7}$ & $9.713\times10^{-8}$ \\
\hline
\end{tabulary}$ $\newline\newline

\begin{tabulary}{\textwidth}{|c|c|c|c|c|c|} 
 \hline
 $j$ & $W_{1}^j$ & $W_{2}^j$ & $W_{3}^j$ & $W_{4}^j$ & $W_{5}^j$\\
 \hline\hline
 $1$ & $4.254\times10^{-5}$ & $1.020\times10^{-7}$ & $2.023\times10^{-7}$ & $1.412\times10^{-7}$ & $2.589\times10^{-7}$ \\
 \hline
 $2$ & $0$ & $2.386\times10^{-7}$ & $2.329\times10^{-7}$ & $8.962\times10^{-8}$ & $1.908\times10^{-7}$ \\
 \hline
 $3$ & $0$ & $2.392\times10^{-7}$ & $1.513\times10^{-7}$ & $4.693\times10^{-8}$ & $1.100\times10^{-7}$  \\
\hline
 $4$ & $0$ & $4.868\times10^{-8}$ & $2.050\times10^{-7}$ & $1.617\times10^{-7}$ & $1.204\times10^{-7}$ \\
\hline
\end{tabulary}

\caption{{\small{Wasserstein distances \eqref{WassersteinBC_SP} between the measured distributions and those predicted in the BC case (top) and SP case (bottom) for the Canneal benchmark in Sec. \ref{subsec:MultiCoreExperiment} at each $\hat{\tau}_{\widehat{\sigma}\in\llbracket 5\rrbracket}$.}}}
\label{tab:WassersteinErrorsBC_SP}
\vspace*{-0.15in}
\end{table}


\section{Conclusions}\label{Sec:Conclusions}
This work explores a new vision for learning and predicting stochastic time-varying computational resource usage from hardware-software profile data in both single and multi-core platforms. We propose to formulate the problem as a distributional learning problem directly in the joint space of correlated computational resources such as processor availability, memory bandwidth, etc. This leads to graph-structured multi-marginal Schr\"{o}dinger bridge problems (MSBPs) where the specific graph structures are induced by underlying single or multi-core nature of computation. At first glance, such formulations for scattered profile data appear computationally intractable because they involve tensorial convex optimization problems where the number of primal variables is exponential in the number of observational snapshots. By leveraging strong duality and graph structures, we show that the computational complexities can be reduced {from exponential to linear without approximation}, and these problems can in fact be solved efficiently via \emph{nonparametric computation} with \emph{maximum likelihood guarantees in the space of distributions}. This enables us to predict the most likely joint (and hence all marginal) computational resource usage distribution at any user-specified query time. We emphasize that our proposed algorithms and the benchmark results reported here, directly work with the scattered profile data without gridding the joint space of computational resource (here, instructions retired, LLC requests, LLC load misses). For the single core case, we illustrate the computational details for a nonlinear model predictive controller benchmark. For the multi-core case, we provide the results for a benchmark software from real-time systems literature.  

The learning-scheduling co-design that builds on the proposed method, will comprise our future work. Another direction of interest is to account for the asynchronous nature of the profiles by formulating the corresponding MSBPs with joint space-time stochasticity, i.e., as distributional generalization of optimal stopping problems initiated in the recent work \cite{eldesoukey2024schr}. Our multi-marginal formulation also opens up the possibility to incorporate additional flow rate constraints as in \cite{dong2024monge}.

\vspace{-1ex}
\appendix

\subsection{Proof of Proposition \ref{PropBCSinkPathCost}}\label{AppProofPropBCSinkPathCost}
With $\bm{K}$ as in \eqref{eqK_BC} and $\bm{U}$ as in \eqref{eqstructuvector}, consider the Hilbert-Schmidt inner product
\begingroup\allowdisplaybreaks
\begin{align}
    &\langle\mathbf K,\mathbf U\rangle = \sum_{\substack{i_{(r,\ell)}, \\(r,\ell)\in\Lambda_{\rm{BC}}}}\left[\bm{K}_{(i_{(r,\ell)} \mid (r,\ell)\in\Lambda_{\rm{BC}})}\right]\left[\bm{U}_{(i_{(r,\ell)} \mid (r,\ell)\in\Lambda_{\rm{BC}})}\right] \nonumber\\
    &= \sum_{\substack{i_{(r,\ell)}, \\(r,\ell)\in\Lambda_{\rm{BC}}}} \Bigg(\bigg(\prod_{\sigma=1}^{s-1}\left[K_{i_{(0,\sigma)},i_{(0,\sigma+1)}}^{0,\sigma}\right]\bigg) \nonumber\\
    &\quad\quad\quad \bigg(\prod_{\sigma=1}^{s}\prod_{j=1}^{J}\left[K_{i_{(j,\sigma)},i_{(0,\sigma)}}^{j,\sigma}\right]\bigg)
    \bigg(\prod_{\sigma=1}^{s}\prod_{j=0}^{J}\left[(\bm{u}_{\sigma}^{j})_{i_{(j,\sigma)}}\right]\bigg)\Bigg) \nonumber\\
    &= \sum_{\substack{i_{(0,\ell)}, \\\ell\in\llbracket s\rrbracket}} \sum_{\substack{i_{(r,\ell)}, \\(r,\ell)\in\llbracket J\rrbracket\times\llbracket s\rrbracket}}\Bigg(\bigg(\prod_{\sigma=1}^{s-1}\left[K_{i_{(0,\sigma)},i_{(0,\sigma+1)}}^{0,\sigma}\right]\bigg) \nonumber\\
    &\quad\quad\quad \prod_{\sigma=1}^{s}\bigg(\prod_{j=1}^{J}\left[K_{i_{(j,\sigma)},i_{(0,\sigma)}}^{j,\sigma}\right]\bigg)
    \bigg(\prod_{j=0}^{J}\left[(\bm{u}_{\sigma}^{j})_{i_{(j,\sigma)}}\right]\bigg)\Bigg) \nonumber\\
    &= \sum_{\substack{i_{(0,\ell)},\:\ell\in\llbracket s\rrbracket}} \Bigg(\bigg(\prod_{\sigma=1}^{s-1}\left[K_{i_{(0,\sigma)},i_{(0,\sigma+1)}}^{0,\sigma}\right]\bigg) \nonumber\\
    &\prod_{\sigma=1}^{s}\sum_{\substack{i_{(r,\ell)}, \\(r,\ell)\in\llbracket J\rrbracket\times\llbracket s\rrbracket}}\bigg(\prod_{j=1}^{J}\left[K_{i_{(j,\sigma)},i_{(0,\sigma)}}^{j,\sigma}\right]\bigg)
    \bigg(\prod_{j=0}^{J}\left[(\bm{u}_{\sigma}^{j})_{i_{(j,\sigma)}}\right]\bigg)\Bigg) \nonumber\\
    &= \sum_{\substack{i_{(0,\ell)},\:\ell\in\llbracket s\rrbracket}} \Bigg(\bigg(\prod_{\sigma=1}^{s-1}\left[K_{i_{(0,\sigma)},i_{(0,\sigma+1)}}^{0,\sigma}\right]\bigg) \nonumber\\
    &\quad\quad\quad\quad\prod_{\sigma=1}^{s}\bigg[\bigg(\underbrace{\bm{u}_{\sigma}^0\odot\Big(\bigodot_{j=1}^{J}K^{j,\sigma}\bm{u}_{\sigma}^{j}\Big)}_{:=\bm{p}_{\sigma}}\bigg)_{i_{(0,\sigma)}}\bigg]\Bigg) \nonumber\\
    &= \bm{p}_1^\top\bigg(\prod_{\sigma=2}^{s-1}K^{0,\sigma-1}\diag(\bm{p}_{\sigma})\bigg)K^{0,s-1}\bm{p}_{s}\label{jmpKUgen} \\
    &= \bm{p}_1^\top\bigg(\prod_{m=2}^{\sigma-1}K^{0,m-1}\diag(\bm{p}_{m})\bigg)K^{0,\sigma-1} \nonumber\\
    &\quad\quad\quad\cdot\diag\bigg(\bm{u}_{\sigma}^0\odot\Big(\bigodot_{j=1}^{J}K^{j,\sigma}\bm{u}_{\sigma}^{j}\Big)\bigg) \nonumber\\
    &\quad\quad\quad\quad\quad\quad\cdot\bigg(\prod_{m=\sigma+1}^{s-1}K^{0,m-1}\diag(\bm{p}_{m})\bigg)K^{0,s-1}\bm{p}_{s} \label{bcKU_marg}
\end{align}
where \eqref{jmpKUgen} follows by the same argument as for Proposition \ref{PropMultiSinkPathCost} (see also \cite[Proposition 1]{bondar2023path}). Then the unimarginal projections \eqref{ProjSimplified_BC_marg0}-\eqref{ProjSimplified_BC_marggen} follow by applying Lemma \ref{lemmaproj1and2} to \eqref{bcKU_marg}.

For the bimarginal projections, let $(j,\sigma)\in\llbracket J\rrbracket\times\llbracket s\rrbracket$ and $\sigma_1,\sigma_2\in\llbracket s\rrbracket$ such that $\sigma_1<\sigma_2$. Starting with \eqref{jmpKUgen}, we have
\begin{align}
    &\langle\mathbf K,\mathbf U\rangle = \bm{p}_1^\top\bigg(\prod_{\sigma=2}^{s-1}K^{0,\sigma-1}\diag(\bm{p}_{\sigma})\bigg)K^{0,s-1}\bm{p}_{s} \nonumber \\
    &= \bm{p}_1^\top\bigg(\bigodot_{m=2}^{s}K^{0,m-1}\bm{p}_m\bigg) \nonumber\\
    &= \bm{1}^\top\diag(\bm{p}_\sigma){K^{0,\sigma-1}}^\top\bigg(\bm{p}_1\odot\bigodot_{\substack{m=2,\\m\neq\sigma}}^{s}K^{0,m-1}\bm{p}_m\bigg) \nonumber\\
    &= \bm{1}^\top\diag\bigg(\bm{u}_{\sigma}^0\odot\Big(\bigodot_{k=1}^{J}K^{k,\sigma}\bm{u}_{\sigma}^{k}\Big)\bigg) \nonumber\\
    &\quad\quad\quad\quad\quad\quad{K^{0,\sigma-1}}^\top\bigg(\bm{p}_1\odot\bigodot_{\substack{m=2,\\m\neq\sigma}}^{s}K^{0,m-1}\bm{p}_m\bigg) \nonumber\\
    &= \bm{1}^\top\diag(\bm{u}_\sigma^0)\diag\bigg(\bigodot_{\substack{k=1 \\ k\neq j}}^{J}K^{k,\sigma}\bm{u}_{\sigma}^{k}\bigg) \nonumber\\
    &\quad\diag\Bigg({K^{0,\sigma-1}}^\top\bigg(\bm{p}_1\odot\bigodot_{\substack{m=2,\\m\neq\sigma}}^{s}K^{0,m-1}\bm{p}_m\bigg)\Bigg)K^{j,\sigma}\diag(\bm{u}_\sigma^j)\bm{1} \nonumber\\
    &= \bm{1}^\top\diag(\bm{u}_\sigma^0) \nonumber\\
    &\diag\Bigg({K^{0,\sigma-1}}^\top\underbrace{\bigg(\bm{p}_1\odot\Big(\bigodot_{\substack{m=2,\\m\neq\sigma}}^{s}K^{0,m-1}\bm{p}_m\Big)\Big(\bigodot_{\substack{k=1 \\ k\neq j}}^{J}K^{k,\sigma}\bm{u}_{\sigma}^{k}\Big)\bigg)}_{=:\bm{\rho}_{(0,\sigma),(j,\sigma)}}\Bigg) \nonumber\\
    &\quad\quad\quad\quad K^{j,\sigma}\diag(\bm{u}_\sigma^j)\bm{1} \label{bcKU2_spect}\\
    &= \bigg(\bm{p}_1^\top K^{0,1}\prod_{k=2}^{\sigma_1-1}\diag(\bm{p}_k)K^{0,k}\bigg) \nonumber\\
    &\quad\cdot\diag(\bm{p}_{\sigma_1})\bigg(\prod_{k=\sigma_1+1}^{\sigma_2-1}K^{0,k-1}\diag(\bm{p}_k)\bigg)K^{0,\sigma_2-1}\diag(\bm{p}_{\sigma_2}) \nonumber\\
    &\quad\cdot\bigg(\prod_{k=\sigma_2+1}^{s-1}K^{0,k-1}\diag(\bm{p}_k)\bigg)K^{0,s-1}\bm{p}_{s} \label{bcKU2_bary}.
\end{align}
Then, \eqref{Proj2Simplified_BC_marg0} follows from \eqref{bcKU2_spect}, and \eqref{Proj2Simplified_BC_margs} follows from expanding $\diag(\bm{p}_{\sigma_1})$ and $\diag(\bm{p}_{\sigma_2})$ in \eqref{bcKU2_bary}, followed by the application of Lemma \ref{lemmaproj1and2}.\qedsymbol
\endgroup


\subsection{Proof of Proposition \ref{PropSPSinkPathCost}}\label{AppProofPropSPSinkPathCost} 
With $\bm{K}$ as in \eqref{eqK_SP} and $\bm{U}$ as in \eqref{eqstructuvector}, we start with the Hilbert-Schmidt inner product
\begingroup\allowdisplaybreaks
\begin{align}
    &\langle\mathbf K,\mathbf U\rangle = \sum_{\substack{i_{(j,\sigma)}, \\(j,\sigma)\in\Lambda_{\rm{SP}}}}\left[\bm{K}_{(i_{(j,\sigma)} \mid (j,\sigma)\in\Lambda_{\rm{SP}})}\right]\left[\bm{U}_{(i_{(j,\sigma)} \mid (j,\sigma)\in\Lambda_{\rm{SP}})}\right] \nonumber\\
    &= \sum_{i_{(j,\sigma)} , (j,\sigma)\in\Lambda_{\rm{SP}}} \left(\prod_{k=1}^{J}\left[\bm{K}_{i_{(1,1)},(i_{(j,\sigma)}\mid (j,\sigma)\in\Lambda_{\rm{SP}}^{k}),i_{(1,s)}}^{k}\right]\right) \nonumber\\
    &\quad\quad\quad\quad(\bm{u}_1^1)_{i_{(1,1)}}\Bigg(\prod_{k=1}^J\prod_{\sigma=2}^{s-1}(\bm{u}_{\sigma}^{k})_{i_{(k,\sigma)}}\Bigg)(\bm{u}_s^1)_{i_{(1,s)}} \nonumber\\
    &= {\bm{u}_1^1}^{\top}\bigodot_{k=1}^{J}\Bigg(\underbrace{K^{k,1}\left(\prod_{\sigma=2}^{s-1}\diag(\bm{u}_{\sigma}^{k})K^{k,\sigma}\right)}_{:=A_k}\Bigg)\bm{u}_s^1 \label{HSForm_SP_1}\\
    &=  {\bm{u}_1^1}^{\top}\Bigg(\bigg( \underbrace{K^{j,1}\bigg(\prod_{m=2}^{\sigma-1}\diag(\bm{u}_{m}^j)K^{j,m}\bigg)}_{:=X_\sigma^j} \nonumber\\
    &\quad\cdot\diag(\bm{u}_{\sigma}^{j})\!\underbrace{K^{j,\sigma}\bigg(\!\prod_{m=\sigma+1}^{s-1}\!\!\!\!\!\diag(\!\bm{u}_{m}^j\!)\!K^{j,m}\!\!\bigg)}_{:=Y_{\sigma}^j}\!\!\bigg)\!\odot\!\underbrace{\bigodot_{k\neq j}A_k}_{:=B_j}\!\Bigg)\!\bm{u}_s^1\nonumber \\
    &=   \tr\Big(\diag({\bm{u}_1^1})\bigr(X_\sigma^j\diag(\bm{u}_\sigma^j)Y_\sigma^j\bigr)\diag(\bm{u}_s^1)B_j^\top\Big) \nonumber\\
    &=   {\bm{1}}^\top\diag\Big(\diag(\bm{u}_\sigma^j)Y_\sigma^j\diag(\bm{u}_s^1)B_j^\top\diag({\bm{u}_1^1})X_\sigma^j\Big) \nonumber\\
    &=  \underbrace{\bm{1}^\top}_{\bm{w}_1^\top}\cdot\diag(\bm{u}_\sigma^j)\cdot\underbrace{\diag\Big(Y_\sigma^j\diag(\bm{u}_s^1)B_j^\top\diag({\bm{u}_1^1})X_\sigma^j\Big)}_{\bm{w}_2}, \label{HSForm_SP_2} 
\end{align}
\endgroup
wherefrom the projections \eqref{ProjSimplified_SP_marg1} and \eqref{ProjSimplified_SP_margs} follow by applying Lemma \ref{lemmaproj1and2} to \eqref{HSForm_SP_1}. Similarly, \eqref{ProjSimplified_SP_marggen} follows by applying Lemma \ref{lemmaproj1and2} to \eqref{HSForm_SP_2}.

Next, for the bimarginal projections, we write
\begingroup\allowdisplaybreaks
\begin{align}
    &\langle\mathbf K,\mathbf U\rangle = {\rm{trace}}\Big( \diag(\bm{u}_1^1)X_{\sigma_1}^j\diag(\bm{u}_{\sigma_1}^{j}) \nonumber\\
    &\quad\quad\quad\quad\quad\quad\quad\quad\quad\cdot Z_{\sigma_1,\sigma_2}^{j}\diag(\bm{u}_{\sigma_2}^j)Y_{\sigma_2}^j\diag(\bm{u}_s^1)B_j^\top\Big)\nonumber \\
    &= \bm{w}_1^\top\diag(\bm{u}_1^1)\underbrace{\Bigg(K^{j,1}\odot\bigg(B_j\diag(\bm{u}_s^1){Y_2^j}^{\top}\bigg)\Bigg)}_{\Phi}\diag(\bm{u}_2^j)\bm{w}_3\label{HS2Form_SP_1}\\
    &= \bm{w}_1^\top\diag(\bm{u}_{s-1}^j)\!\underbrace{\Bigg(\!\!K^{j,s-1}\!\odot\!\bigg(\!\!{X_{s-1}^j}^{\!\!\!\!\!\!\top}\diag(\bm{u}_1^1)B_j\!\!\bigg)\!\!\Bigg)}_{\Phi}\!\diag(\bm{u}_s^1)\bm{w}_3\label{HS2Form_SP_s}\\
    &= \bm{w}_1^\top\diag(\bm{u}_{\sigma_1}^j)\underbrace{\Bigg(\!\!Z_{\sigma_1,\sigma_2}^j\!\odot\!\bigg(\!\!{X_{\sigma_1}^j}^{\!\!\top}\!\diag(\bm{u}_1^1)B_j\diag(\bm{u}_s^1){Y_{\sigma_2}^j}^{\!\!\top}\!\bigg)\!\!\Bigg)}_{\Phi} \nonumber\\
    &\quad\quad\quad\quad\quad\quad\cdot\diag(\bm{u}_{\sigma_2}^j)\bm{w}_3,\label{HS2Form_SP_gen}
\end{align}
\endgroup
where $\bm{w}_1=\bm{w}_3=\bm{1}$. Now \eqref{Proj2Simplified_SP_marg1} follows from \eqref{HS2Form_SP_1}, \eqref{Proj2Simplified_SP_margs} from \eqref{HS2Form_SP_s}, and \eqref{Proj2Simplified_SP_marggen} from \eqref{HS2Form_SP_gen}.\qedsymbol

\bibliographystyle{IEEEtran} 
\bibliography{TCST-MSBP-CPSfrontier.bib}

\end{document}